\documentclass{article}
\usepackage{amsmath}
\usepackage{theorem}
\usepackage{amssymb}
  \usepackage[dvips]{graphics}

\newtheorem{Proposition}{Proposition}
  \newtheorem{Remark}{Remark}
  \newtheorem{Corollary}[Proposition]{Corollary}
  \newtheorem{Lemma}[Proposition]{Lemma}
  
  \newtheorem{Theorem}[Proposition]{Theorem}

\newtheorem{Definition}[Proposition]{Definition}

\newlength{\strikewidth}
\newlength{\strikelength}
\newenvironment{proof}{\par\noindent {\sc Proof.}}{$\Box$}

\newcommand{\be}{\begin{equation}}
\newcommand{\ee}{\end{equation}}

\def\bfY{\mathbf{Y}}
\def\bfm{\mathbf{m}}
\def\bfj{\mathbf{j}}
\def\bfd{\mathbf{d}}
\def\bfT{\mathbf{T}}
\def\bor{\mathcal{B}}
\def\bfy{\mathbf{y}}
\def\bfk{\mathbf{k}}
\def\z{\noindent}
\def\pint{{- \kern-1.1em\int}}
\def\pist#1#2{\noindent\hangindent 2em\hangafter1\hbox to 2em{#1\hfil~~}#2}
\def\cal{\mathcal}
\def\e{\mathrm{e}}

\def\Box{{\hfill\hbox{\enspace${\sqre}$}} \smallskip}
\def\sqr#1#2{{\vcenter{\vbox{\hrule height .#2pt
                             \hbox{\vrule width .#2pt height#1pt \kern#1pt
                                   \vrule width .#2pt}
                             \hrule height .#2pt}}}}
\def\sqre{\mathchoice\sqr54\sqr54\sqr{4.1}3\sqr{3.5}3}

\def\CC{\mathbb{C}}
\def\RR{\mathbb{R}}
\def\ZZ{\mathbb{Z}}
\def\NN{\mathbb{N}}
\def\erm{\mathrm{e}}
\def\lap{\mathcal{L}}
\begin{document}

\title{Complex Singularity Analysis for a Nonlinear PDE }
\author{{\it O. Costin}\\ \small{Math Department, Rutgers University}\\
{\it S. Tanveer} \\
\small{Math Department, The Ohio State University}}

\maketitle
\date{}
\bigskip

\begin{abstract}
  We introduce a method of rigorous analysis of the location and type of
  complex singularities for nonlinear higher order PDEs as a function of the
  initial data.  The method is applied to determine rigorously the asymptotic
  structure of singularities of the modified Harry-Dym equation
  $$ H_t + H_y = - \frac{1}{2} H^3 + H^3 H_{yyy} ~~: ~~H(y, 0) = y^{-1/2} $$
  for small time at the boundaries of  the  sector of analyticity. Previous work
  \cite{CPAM}, \cite{invent03} shows existence, uniqueness and Borel
  summability of solutions  of general PDEs. It is shown that the 
  solution  to the above initial value problem is represented
  convergently by a series in a fractional power of $t$ down to a small
  annular neighborhood of a singularity of the leading order equation. We
  deduce that the exact solution has a singularity nearby  having, to
  leading order, the same type.

\end{abstract}

\section{Introduction} 

The singularity structure of solutions of nonlinear partial differential
equations in the complex plane is not a well understood subject.
Insofar as the authors are aware, except for integrable cases, there are
no general methods in the literature to analyze the type and location of
singularities of solutions of nonlinear systems of PDEs in terms of the
initial data.

The goal of the present paper is to develop a relatively general and
constructive technique to address this issue, which applies to equations that
admit formal asymptotic solutions.

In view of the complexity of the analysis, and for more clarity, we describe
and apply the method on a number of concrete problems. It will be however
transparent that the method is much more general.

Formal asymptotic similarity solutions could exist for small or large time, or
when one approaches the finite blow-up time of a similarity solution of a PDE.
We prove that complex singularities of these formal asymptotic solution
actually correspond to singularities of the solution of the full PDE.

The motivation for understanding complex singularity formation of PDEs, aside
from intrinsic interest, is that in some cases of physical interest
\cite{Moore}, \cite{Tanveer93b}, \cite{Cowley},  there is evidence  
that singularities that appear in the real physical domain after a finite time
can be traced to the complex plane.

The procedure consists in the following steps: (i) an early time asymptotic
expansion in powers of $t$, the validity of which is justified for the
modified Harry-Dym equation in \cite{invent03}, (ii) introduction of
appropriately scaled ``inner" dependent and independent variables beyond the
region of validity of the expansion (i), (iii) determination of singularities
of the leading order equation and (iv) proof that a secondary expansion in
scaled time, involving inner-variables, is convergent in a domain encircling a
singularity of the leading order solution.  Insofar as analysis of the leading
order equation (in step (iii) above) is concerned, which (typically a
nonlinear ODE), formal calculations have been used before (see
\cite{Tanveer93}, and references in \cite{Invent}).  These can now be
rigorously derived from the general theory introduced in \cite{Invent}.

The present paper justifies the above four step procedure for the
modified Harry-Dym equation This equation arises in the small surface tension
limit of Hele-Shaw interfacial evolution \cite{Tanveer93} in the neighborhood
of an initial zero of the derivative of  an associated conformal
map. The justification of singularity formation is a crucial first-step to
understanding ``daughter"-singularity phenomena where a smoothly evolving
interface corresponding to a zero-surface tension solution is singularly
perturbed in $O(1)$ time by arbitrarily small surface tension effect.

Consider the following initial value problem 
for the modified Harry-Dym equation:
\begin{equation}
\label{0.2.1}
\frac{\partial H}{\partial t }
+ \frac{\partial H}{\partial y} 
- H^3 \frac{\partial^3 H}{\partial y^3} - \frac{H^3}{2} = 0
~~{\rm with}~~H(y, 0) = y^{-1/2}
\end{equation}
 Theorem 36 and Corollary 37 in \cite{invent03} imply that  
for any $t \in [0, T]$, for large
enough $|y-t|/t^{2/9}$ with 
$\arg (y-t) \in \left ( - \frac{4}{9} \pi, \frac{4}{9} \pi \right )$,
there exists a unique solution $H(y, t)$ to (\ref{0.2.1})
with $H(y, t) \sim y^{-1/2}$ as $|y-t|/t^{2/9} \rightarrow \infty$,
with the following asymptotic expansion
for $t \ll 1 $: 
\begin{equation}
\label{0.2.3}
H(y,t) = (y-t)^{-1/2} \sum_{n=0}^\infty P_n \left ( 
\frac{t}{(y-t)^{9/2}}, 
\frac{t}{(y-t)} \right )  
\end{equation}
where $P_0 =1$ and $P_n $  is a homogeneous polynomial determined
recursively in terms of $P_{n-1}$, $P_{n-2}$,...$P_1$.  The first two  polynomials are
\begin{equation}
P_1 (a, b) = -\frac{15}{8} a - \frac{1}{2} b  ,~~~~\\ 
P_2 (a, b) = \frac{25875}{128} a^2 + \frac{195}{32} a b + \frac{3}{8} b^2
\end{equation}
Further, if we introduce  the  scaled
variables 
\begin{equation}
\label{eq:etaG}
\eta = \frac{x-t}{t^{2/9}} ~~;~~ \tau = t^{7/9} ~~;~~ 
H(y (\eta, t), t) = t^{1/9} G(\eta, \tau), 
\end{equation}
then, according to 
Corollary 37 in \cite{invent03}\footnote{The variable
variable $\zeta= \eta^{3/2}$ is used there.},
for $|\eta|$ sufficiently large, with 
$\arg \eta \in \left ( -\frac{4}{9} \pi, ~\frac{4}{9} \pi \right ) $,
the function $G$  has a convergent series expansion in $\tau$: 
\begin{equation}
\label{0.2.6}
G(\eta, \tau) = \sum_{k=0}^\infty \tau^k G_k (\eta) 
\end{equation}

\vskip 0.2truein

In this paper, among other results, it will be
shown that the convergence of the series (\ref{0.2.6})
actually holds in
an extended domain in $\eta$ that includes
at least a region close to a singularity ${\hat \eta}_s$ of $G_0 (\eta)$
in a neighborhood of the boundary  $\arg \eta \in \left (-\frac{4}{9} \pi, 
\frac{4}{9} \pi \right )$ for large $|{\hat \eta}_s|$.  

Substituting (\ref{eq:etaG}) into
(\ref{0.2.1}), we obtain the following equation for 
$G(\eta, \tau)$:
\begin{equation}
\label{0.2.5}
-\frac{G}{9} - \frac{2}{9} \eta G_\eta + \frac{7}{9} \tau G_\tau + 
\frac{\tau}{2} G^3 - G^3 G_{\eta \eta \eta} =0
\end{equation}
From (\ref{0.2.6}), it follows that 
the leading order solution $G_0$ satisfies
\begin{equation}
\label{0.2.7}
\frac{1}{9} G_0 
+ \frac{2}{9} \eta G_0^\prime  + G_0^3 G_0^{\prime \prime \prime} = 0 
\end{equation}
In order for $G$ in  
(\ref{0.2.6}) to match the asymptotic expansion (\ref{0.2.3})  we need to
require that 
\begin{equation}
\label{0.2.11}
G_0 (\eta)=\eta^{-1/2}(1+o(1));\ \  |\eta|\ \text{large,}~\arg \eta \in \left ( -\frac{4 \pi}{9}, \frac{4 \pi}{9} \right )  
\end{equation}
The solution $G_0 (\eta)$ to the
leading order ODE (\ref{0.2.7}) with asymptotic condition (\ref{0.2.11}) 
have been studied before.
Numerical solutions were found \cite{Tanveer93} and
computational evidence suggested that there is a cluster of singularities
${\hat \eta}_s$, where $G_0 (\eta)\sim 
~\e^{i \pi/3} \left (\frac{\eta_s}{3} \right )^{1/3} (\eta
- {\hat \eta_s} )^{2/3} $.  
Using the fact that $G_0(\eta)$ is indeed a similarity
solution to the Harry-Dym equation, which is integrable, it was shown
\cite{Fokas} that (\ref{0.2.7}) can be transformed    
to  Painlev\'e P$_{\rm II}$. 
Isomonodromic methods were used to
prove existence and uniqueness of sectorially analytic solution
for $\arg \eta \in \left ( -\frac{4}{9} \pi, \frac{4}{9} \pi \right )$ 
that satisfies far-field condition 
(\ref{0.2.11}).
Outside this sector, the
behavior of the transformed equation solution 
is given by elliptic functions,
whose singularities 
can be related to the two-thirds singularity of $G_0 (\eta)$, as above.

However, unlike the isomonodromic method, 
the method based on
generalized Borel summation summation \cite{IMRN}, \cite{DMJ} 
applies to
initially small solutions of non-integrable equations as well. 
To apply this analysis
in our example, which does not satisfy all the conditions in
\cite{Invent}, small adaptations of the proofs are needed (see
Appendix).  One can determine that for large $\eta$, uniformly in the
sector $\arg \eta ~\in ~\left [- \frac{4 \pi}{9} - \delta, ~\frac{4
    \pi}{9} - \delta \right ]$ for some $ \delta ~\in \left (0,
  \frac{2}{9} \pi \right ) $, except for an exponentially small region
around singularity of $G_0$, the asymptotic series of $G_0 (\eta)$ is
of the form
\begin{equation}
\label{0.3.1}
G_0 (\eta) ~\sim~ \eta^{-1/2} 
U (\zeta )   + O(\eta^{-5})
\end{equation}
where 
\begin{equation}
\label{0.3.1.5}
\zeta = -\log ~C +\frac{9}{8}~\log~ \eta +  \frac{i 4\sqrt{2}}{27} \eta^{9/4}+(2 {\hat n} -1) i \pi
\end{equation}
with the principal branch of the log, where $C$ is a  Stokes
constant of $G_0 $ in the large $\eta$ expansion for $\arg \eta \in \left [ - \frac{4}{9} \pi + \delta,
\frac{4}{9} \pi - \delta \right ]$
\begin{equation}
\label{0.3.1.6}
G_0 (\eta) \sim \eta^{-1/2} 
\left [ 1 + \sum_{k=1}^\infty a_k \eta^{-9 k/2} \right ]
\end{equation}
The function $U(\zeta) $ is determined implicitly from the equation:
\begin{equation}
\label{0.3.2}
\zeta = \log 4 - 2 - i \pi - 2 \sqrt{U} - \ln \left ( \frac{1-\sqrt{U}}{1+\sqrt{U}} \right )
\end{equation} 
and $ U(\zeta)$ has a singularity (where $U=0$) at 
$\zeta = \zeta_s \equiv \log 4 - 2 -i \pi $, corresponding to a string of
singularities at $\eta=\eta_s$, where 
\begin{equation}
\label{0.3.2.1}
\frac{i 4\sqrt{2}}{27} \eta_s^{9/4} +\frac{9}{8} \log \eta_s = -2 + 
\log 4 - 2 {\hat n} i \pi + \log C 
\end{equation}
where ${\hat n}\in\NN$ has to be large for $\eta_s$ to be large.  For
large $|\eta_s|$(large ${\hat n}$), it is to be noted that $\arg
\eta_s$ is close to $-\frac{4 \pi}{9}$, the anti-Stokes line. There is
similarly another quasi-periodic array of singularities close to $\arg
\eta = \frac{4 \pi}{9}$, but our focus will be only on the ones in the
lower-half plane. It can be shown that for large $|\eta_s|$ the
singularities of $G_0$ lie within an exponentially small distance of
$\eta_s$ and, to leading order, are of the same type. This can be
further verified directly from the equation for $G_0$.

\begin{Remark}\label{R1}
It is easy to check that $G_0$ cannot be zero, except at a singularity
$\eta={\hat \eta}_s$. Furthermore , in any domain $\mathcal{D}$ that
excludes a neighborhood of the singularities of $G_0$, and extends to $\infty$
so that $\arg \eta \in \left [ -\frac{4}{9} \pi + \delta, \frac{4}{9}
\pi + \delta \right ] $, it follows from differentiability of the 
asymptotics of solutions of ODEs \cite{Wasow} that
\begin{equation}
\label{eqderG0}
\sup_{\eta \in \mathcal{D}} |\eta^{1/2} G_0 (\eta)| ,~~\sup_{\eta
  \in \mathcal{D}} |\eta^{7/2} G_0^{\prime \prime \prime} (\eta) |
<C 
\end{equation}
\end{Remark}
 
\begin{Remark}\label{U01}

By (\ref{0.3.2}), near the singularity $\eta=\eta_s$ we have
\begin{equation}
\label{2.2.6.7}
\frac{2}{3} U^{3/2} + O(U^{5/2}) = \zeta - \zeta_s ~=~\frac{1}{\eta_s} \left \{ \frac{9}{8} 
+ \frac{i \sqrt{2}}{3} \eta_s^{9/4} \right \} (\eta-\eta_s) 
\left [ 1 + O \left (\eta_s^{-1} (\eta-\eta_s) \right ) \right ]  
\end{equation}
and hence for $\eta-\eta_s = o(\eta_s^{-5/4})$, for large enough $|\eta_s|$, 
\begin{equation}
\label{2.2.6.8}
U ~\sim ~\e^{i \pi/3} \left ( \frac{\eta_s^{5/6}}{2^{1/3}} \right ) 
\left ( 1 - \frac{27 i}{8 \sqrt{2}} \eta_s^{-9/4} \right )^{2/3} (\eta-\eta_s)^{2/3}  
\end{equation}
Note that if $r_i |\eta_s|^{-5/4} < |\eta-\eta_s| < r_0 |\eta_s|^{-5/4}$, with
$r_0 > r_i$ small, then there exists upper and lower bounds for $|U|$,
independent of $\eta_s$ for large $|\eta_s|$.  Since the singularity ${\hat
\eta}_s$ of $G_0$ is exponentially close to $\eta_s$, it follows that
the lower bound of $G_0$ in this annular region is also independent of
$|\eta_s|$.

\end{Remark}

Given these leading order singularities for $G_0(\eta)$,  we 
investigate the series expansion (\ref{0.2.6}),  known 
to converge for large enough $|\eta|$ in any compact subset of $\arg
\eta \in \left ( -\frac{4}{9} \pi, \frac{4}{9} \pi \right )$, in a
neighborhood of a singularity of $G_0$. 
 \begin{Remark}
  The domain $\mathcal{D}$ in the next theorem, sketched in Fig. 1, is
  specified in Definition \ref{defD_A}. It contains a small
  annular region of a singularity $\eta_s$ of $U$ (cf.
  (\ref{0.3.2.1}), (\ref{0.3.2})) for large $\hat{n}\in\NN$ and a
  sectorial region $\arg \eta \in \left ( -\frac{2\pi}{9} + \delta ,
    \frac{2\pi}{9} - \delta \right )$, for $\delta$ small and  large
  $\eta$.
\end{Remark} 
\z Comparing powers of $\tau^k$ (for $k \ge 1$) obtained
by substituting power series (\ref{0.2.6}) into (\ref{0.2.5}) results in:
\begin{equation}
\label{0.2.8} 
G_0^3 {\mathcal L}_k G_k = R_k 
\end{equation}
where the linear operator $\mathcal{L}_k$ is defined by
\begin{equation} 
\label{0.2.9}
\mathcal{L}_k u = u^{\prime\prime\prime} + \frac{2}{9 G_0^3} \eta u^\prime - \left (\frac{\beta_k}{G_0^3}
+ \frac{3 G_0^{\prime \prime \prime}}{G_0} \right ) u  ~~{\rm where}~~\beta_k = \frac{7 k -1}{9}
\end{equation}
and the terms $R_k $ on the right side of (\ref{0.2.8}) are given by
\begin{equation}
\label{0.2.10}
R_k (\eta) = \frac{1}{2} \sum_{\sum k_i= k-1} G_{k_1} G_{k_2}
G_{k_3} - \sum_{k_j < k, \sum k_j = k} G_{k_1} G_{k_2} G_{k_3} G_{k_4}^{\prime \prime \prime} 
\end{equation}
In order to match to the asymptotic expansion 
expansion (\ref{0.2.3}), we  require
\begin{equation}
\label{0.2.12}
G_k (\eta) ~\sim ~\frac{A_k}{\eta^{k+1/2}} ;\ \  |\eta|\ \text{large,}~\arg \eta \in \left ( -\frac{4 \pi}{9}, \frac{4 \pi}{9} \right )  
\end{equation}
for some specific constants $A_k$ ($A_1 = -1/2$,
$A_2 = \frac{3}{8}$, $A_3 = -\frac{5}{16}$,...).  As explained later, it
is not necessary to impose
(\ref{0.2.12}); any solution $G_k$ which 
approaches 0 as $|\eta| \rightarrow \infty$ with
$\arg \eta \in \left ( - \frac{2}{9} \pi - \delta, \frac{2}{9} \pi 
+ \delta \right ) $ at a rate faster than $\eta^{-1/2}$ must necessarily
have the asymptotic behavior
(\ref{0.2.12}) (See Remark \ref{R4}).

\begin{Theorem}
\label{T0.1}
The expansion (\ref{0.2.6}) is convergent in $\mathcal{D}$ for all
sufficiently small $\tau$. In particular, for any singularity $\hat \eta_s$ of
$G_0 (\eta)$ near the anti-Stokes line $\arg \eta = -\frac{4}{9}\pi$ with
$|\hat \eta_s|$ sufficiently large, there is a singularity of $ G(\eta, \tau)$
for small $\tau$, to leading order of the same type, approaching it as $\tau
\rightarrow ~0^+$ .

\end{Theorem}

\begin{Remark}
  The convergence of the Taylor expansion in $\tau$ and the bounds on
  $G_k$ and $G_k^\prime$ suffice to show that $G(\eta, \tau)$ has the
  singularities close to those of $G_0 (\eta)$ since for a circle
  $S_{\epsilon_1}$ of radius $\epsilon_1$ around $\eta_s$ we have
\begin{equation}
\label{0.3.3}
\frac{1}{2 \pi i} \oint_{S_{\epsilon_1}} \frac{G_\eta}{G} d \eta ~\sim ~ 
\frac{1}{2\pi i} \oint_{S_{\epsilon_1}} 
\frac{G_0^\prime}{G_0} d \eta ~\sim~\frac{2}{3} + 
O \left (\epsilon_1^{1/3} , \tau \right ) 
\end{equation}
 For small $\tau$, $G (\eta, \tau)$  thus  has, to leading order,
  a branch-point of algebraic order $2/3$. 
\end{Remark}

\begin{Remark}
 The convergence of the series (\ref{0.2.6}) in $\mathcal{D}$ is a corollary of the following lemma.
\end{Remark}

\begin{Lemma}
\label{L0.2}
There exist constants $A$ and $B$ independent of $j \ge 1$,  with $A > 1$, 
$ 0 < B < 1$,  so that 
\begin{equation}
\label{0.3.4}
\| \eta^{3/2} G_j \|_{\infty, \mathcal{D}} ~\le~\frac{B A^j}{j^3}
\end{equation}
\begin{equation}
\label{0.3.5}
\| \eta^{5/2} G_j^\prime \|_{\infty, \mathcal{D}} ~\le~\frac{B A^j}{j^2}
\end{equation}
\begin{equation}
\label{0.3.6}
\| G_j^{\prime\prime\prime} \|_{\infty,\mathcal{D}} ~\le~\frac{B A^j}{j^2}
\end{equation}
\end{Lemma}

\begin{Remark} 
\label{R4}
The proof of this key Lemma that leads to the proof of Theorem
\ref{T0.1} is given at the end of \S6.  First, we prove a Lemma
bounding the $R_k (\eta)$. This provides bounds of $G_k$ using a
suitable inversion of $\mathcal{L}_k$ in (\ref{0.2.9}).  The estimates
suffice for our purpose but are not sharp, as (\ref{0.2.12}) implies a
faster decay rate in $\eta$.  
The uniqueness of the solution $G(\eta, \tau)$ in
the regime $|\eta| \gg 1$ for $\arg \eta \in \left ( -\frac{4}{9} \pi,
  \frac{4}{9} \pi \right )$ with $G (\eta, \tau) \sim \eta^{-1/2} $ is
shown in \cite{invent03}.  
\end{Remark}
The proof of Lemma~\ref{L0.2} is by induction; we first prove an 
general result for sums of type (\ref{0.2.10}).
\begin{Lemma}
 \label{L1}
With $G_0$ the same as before, there is a constant $K_3$ so that for any $A \in
(1,\infty)$, $B \in (0, 1)$, $\NN\ni k\ge 2$
and $\{G_j\}_{j=1,...,k-1}$ defined in $\mathcal{D}$ 
which satisfy (\ref{0.3.4})-(\ref{0.3.6}) we have
in (\ref{0.2.10}),
\begin{equation}
\label{0.3.7}  
\| \eta^{3/2} R_k \|_{\infty, \mathcal{D}} 
~\le ~\frac{K_3}{k^2}~ (B^2 A^k + B A^{k-1})
\end{equation}
\end{Lemma}

\begin{proof}
It is convenient to break up $R_k$  as:
$$ R_k = R_{0,k} + R_{1,k} $$
where for $k=1$,
$$R_{0,k} = \frac{G_0^3}{2} $$
and for $k > 1$,
\begin{multline*}
\label{2.2.9.1}
%\begin{eqalign}
R_{0,k} = \frac{3}{2} G_0 \sum_{*,k-1} G_{k_1}
G_{k_2} + \frac{3}{2} G_0^2 G_{k-1} 
- 3 G_0 G_0^{\prime \prime \prime} \sum_{*,k}
G_{k_1} G_{k_2} 
- G_0^{\prime \prime \prime} \sum_{*,k}
G_{k_1} G_{k_2} G_{k_3} 
\\
-3 G_0^2 \sum_{*,k} G_{k_1} G_{k_2}^{\prime\prime\prime}
- 3 G_0 \sum_{*,k} 
G_{k_1} G_{k_2} G_{k_3}^{\prime \prime \prime}
%\end{eqalign}
\end{multline*}
where $\displaystyle \sum_{*,\alpha}$ denotes summation over
$k_i\ge 1$ with $\sum_i k_i=\alpha $ and
$$ R_{1,k} (\eta) = \frac{1}{2} \sum_{*,k-1} G_{k_1} G_{k_2}
G_{k_3} - \sum_{*,k} G_{k_1} G_{k_2} G_{k_3} G_{k_4}^{\prime \prime \prime} $$

The proof follows by using the upper bounds on $G_j$, $G_j^\prime$ and
$G_j^{\prime\prime\prime}$ in (\ref{0.3.4})-(\ref{0.3.6}) for 
$k-1 \ge j \ge 1$,  using (\ref{eqderG0})
and noting that 
$$  
\sup_{k}\left\{\sum_{*,k-1} \frac{k^3}{k_1^3 k_2^3},~ \sum_{*,k} \frac{k^2}{k_1^3 k_2^2},~
 \sum_{*,k} \frac{k^2}{k_1^3 k_2^3 k_3^2},~
\sum_{*,k} \frac{k^3}{k_1^3 k_2^3 k_2^3 },~
 \sum_{*,k} \frac{k^2}{k_1^3 k_2^3 k_3^3 k_4^2 }\right\}<\infty
$$

\end{proof}
\section{Proofs}
The proofs rely on bounding $G_k$ in (\ref{0.2.6}).  For given
$k_0$ and $1 \le k \le k_0$, it can be seen that the 
solution $G_k$ to $\mathcal{L}_k G_k = \frac{R_k}{G_0^3}$ that goes to
0 as $\eta\to\infty$ in the sector $\arg \eta \in \left (-\frac{2}{9}
  \pi + \delta , \frac{2}{9} \pi - \delta \right ) $, with $0 < \delta
< \frac{\pi}{63} $, is given by
\begin{equation}
\label{1.4.13.0}
G_k (\eta) = \sum_{j=1}^3 u_j (\eta) 
\int_{\infty \e^{\theta_j}}^\eta v_j (\eta') \frac{R_k (\eta')}{G_0^3 (\eta')} d \eta' 
\end{equation}
Here $\theta_1 = -\frac{2}{9} \pi +
\delta $, $\theta_2 = \frac{2}{9} \pi -\delta $ and $\theta_3 = 0$.
Also, in (\ref{1.4.13.0}), $u_1$, $u_2$ and $u_3$
are three independent solutions of $\mathcal{L}_k u = 0$, with the following
asymptotic behavior for large $\eta$ (see \cite{Wasow}):
$$ u_1 (\eta) \sim ~\eta^{-15/8} \exp \left [i \frac{4 \sqrt{2}}{27}  \eta^{9/4} \right ] $$ 
$$ u_2 (\eta) \sim ~\eta^{-15/8} \exp \left [-i \frac{4 \sqrt{2}}{27}  \eta^{9/4} \right ] $$
$$ u_3 (\eta) \sim ~\eta^{\frac{9}{2} \beta_k } $$
 where $\beta_k$ is defined in (\ref{0.2.9}), 
$(v_1, v_2, v_3)^T$ is the third column of $\Phi^{-1}$ and
$$ \Phi (\eta)  = \left [ \begin {array}{ccc} u_1 & u_2 & u_3 \\
           u_1^\prime & u_2^\prime & u_3^\prime \\
           u_1^{\prime \prime} & u_2^{\prime \prime}  & u_3^{\prime \prime} 
           \end {array} \right ]
$$
It is easily seen that for large $|\eta|$ we have
$$ v_1 (\eta) \sim ~-\frac{9}{4}~\eta^{-5/8} \exp \left [-i \frac{4 \sqrt{2}}{27} \eta^{9/4} \right ] $$ 
$$ v_2 (\eta) \sim ~-\frac{9}{4}~\eta^{-5/8} \exp \left [i \frac{4 \sqrt{2}}{27}  \eta^{9/4} \right ] $$
$$
v_3 (\eta) \sim ~\frac{9}{2} \eta^{-\frac{9}{2} \beta_k -\frac{5}{2}} $$
and
that $v_1$, $v_2$, $v_3$ are three independent solutions to the adjoint third
order linear equation $\mathcal{L}_k^{+} v = 0$, where the coefficients are
regular when $G_0 \ne 0$.  The $G_k$ determined
from (\ref{1.4.13.0}) are bounded in any domain that excludes the singularities
of $G_0$ (the only places where $G_0 = 0$), and it is easily seen
that the bounds (\ref{0.3.4})-(\ref{0.3.6}) are valid for $1 \le k \le
k_0$ for large $A$ and $B$ (depending
on $k_0$). 

To prove the bounds (\ref{0.3.4})-(\ref{0.3.6}) in Lemma \ref{L0.2}
for all $k$, it is sufficient to prove them for sufficiently large $k
$ (large $\beta_k$).  

\z {\bf Note}. We have to treat separately two regimes of $\eta$ and
$\beta_k$ : (a) $|\eta| =O(\epsilon \beta_k^{4/9})$ or larger and (b)
$|\eta| = o(\beta_k^{4/9})$. These require different integral
representation of $G_k$ and choice of domain.

\subsection{Control in the regime (a), 
$|\eta| >~{\rm constant}~ \beta_k^{4/9}$} 

It is convenient to define 
$$ \chi = \beta_k^{-4/9} \eta ~{\rm and} ~z_k (\chi) = G_k (\beta_k^{4/9} \chi) $$
Then, using variation of parameters (see \S5), we have 
$$ z_k (\chi) = \mathcal{\tilde V} \left [ {\hat R} \right ] (\chi) $$
where
\begin{equation}
\label{1.4.13}
\mathcal{\tilde V} [\hat R] (\chi) \equiv 
\sum_{j=1}^3 \frac{1}{\beta_k^2} \int_{\infty \e^{i \theta_j}}^\chi  
\e^{\beta_k ~\left [P_j (\chi) - P_j ({\tilde \chi}) \right ] + W_j (\chi)
~-~W_j ({\tilde \chi}) } n_{j,3} ({\tilde \chi}) {\hat R} ({\tilde \chi}) d{\tilde \chi}    
\end{equation}
where ${\hat R}$ depends on $R_k $ and $z_k$; $n_{j,3}$,   and  $W_j$ are given
functions of $\chi$, whose exact expression is irrelevant, with behavior
$n_{j,3} = O(\chi^{-5/2})$, $W_1 ~=~-\frac{15}{8} ~\ln ~\chi + o(1)$, $W_2
~=~-\frac{15}{8} ~\ln ~\chi + o(1)$ and $W_3 = o(1)$ for large $\chi$; and $
P_1^\prime (\chi) $, $P_2^\prime (\chi) $ and $P_3^\prime (\chi)$ are the
three roots of the cubic
\begin{equation}
\label{1.4.13.3}
\alpha^3 + \frac{2}{9} \chi^{5/2} \alpha - \chi^{3/2} = 0 
\end{equation}
with the following asymptotic behavior for large $\chi$:
$$P_1 ~= ~\frac{4 \sqrt{2}}{27} i \chi^{9/4} - \frac{9}{4} \ln ~\chi + o(1) $$
$$P_2 ~=~-\frac{4 \sqrt{2}}{27} i \chi^{9/4} - \frac{9}{4} \ln ~\chi + o(1) $$
$$ P_3 ~=~\frac{9}{2} \ln ~\chi + o(1) $$
as $\chi \rightarrow \infty$ for 
$\arg \chi \in \left ( -\frac{2}{9} \pi + \delta , \frac{2}{9} \pi - \delta \right )$,
with $0< \delta < \frac{\pi}{63} $. 
In (\ref{1.4.13}), $\theta_1 = -\frac{2}{9} \pi + \delta $, 
$\theta_2 = \frac{2}{9} \pi - \delta $ and $\theta_3 = 0 $. 

It is necessary that the operators $\mathcal{\tilde V}$ be defined in
a suitable domain $\mathcal{E}$ in the $\chi$-plane 
  containing the integration path where the bounds for the
previous $z_j$, $j=1,..., k-1$ are available to  estimate  $R_k$.
Also $\mathcal{\tilde V}$ need to be bounded for large $\beta_k$.  To
satisfy the latter  requirement for each $j$ ($j=1,2,3$), any point
$\chi \in \mathcal{E}$ must have the property that it can be connected
to $ \infty \e^{i \theta_j}$  along  a path $\mathcal{{\tilde C}}_j$
entirely in $\mathcal{E}$ so that on the path ${\tilde \chi} (s)$,
parameterized by the arclength $s$ increasing towards $\infty$,
$$ \frac{d}{ds} ~\Re P_j ({\tilde \chi} (s)) \ge ~C |{\tilde \chi} (s) |^{5/4} ~>~0 ~,~~{\rm for} ~~
j = 1, 2 $$
$$
\frac{d}{ds} ~\Re P_3 ({\tilde \chi}(s)) \ge ~C~|{\tilde \chi} (s) |^{-1} ~>~0~,~ $$
where
$C$ is a constant independent of $\chi$. It is shown in \S3 that these
properties are ensured if we choose
$$ \mathcal{E} = \left \{ \chi :\,\,\chi ~{\rm to~the~right~of}~{\partial
\mathcal{E}_L },~\arg \chi \in \left (-\frac{2}{9}\pi + \delta, ~\frac{2}{9}
\pi - \delta \right ) \right \} $$ \z  where  $\partial \mathcal{E}_L $
is the polygonal line connecting $\chi_1$, $\chi_3$ and $\chi_2$, 
and  where
$$ \chi_3 = \epsilon ~~,~~\chi_2 = \chi_3 + {\tilde \rho} \e^{i 2 \pi /3}
~~,~~\chi_1 = \chi_3 + {\tilde \rho} \e^{-i 2 \pi/3} $$  Here  ${\tilde
\rho}$ is chosen so that $\arg \chi_1 = -\frac{2 \pi}{9} + \delta $ and $\arg
\chi_2 = \frac{2}{9} \pi - \delta$ and $\epsilon$  is suitably small,
independent of $k$, so that ${\tilde \delta}$ appearing in the proof of
Theorem \ref{T5.1} in \S 6 is smaller than $\frac{1}{2}$.  The domain
$\mathcal{E}$ is sketched in Fig. 3.  Corresponding to $\mathcal{E}$, we
define the domain $\mathcal{E}_k$ (Figure 2)
$$ \mathcal{E}_k = 
\left \{ \eta ~:~ \beta_k^{-4/9} \eta = \chi \in \mathcal{E} \right \} $$

\subsection{Control in regime (b), $\eta =o\left(\beta_k^{4/9}\right)$:}\label{S:WKB}

In this case, as shown in \S6, we can write
$$ G_k (\eta) = \mathcal{V} \left [ {\hat R}_k \right ] + \sum_{j=1}^3 a_j g_j (\eta) $$
where 
\begin{equation}
\label{1.5}
\mathcal{V} \left [{\hat R}_k \right ] \equiv
\sum_{j=1}^3 \frac{\beta_k^{-2/3}}{3} \omega_j G_0 (\eta) 
\int_{\eta_{j,k}}^\eta G_0 (\eta') {\hat R}_k (\eta') 
\e^{\omega_j \beta_k^{1/3} [P (\eta)- P(\eta') ]} d\eta' 
\end{equation}
$g_j = G_0 ~\e^{\omega_j \beta_k^{1/3} P} $,
${\hat R}_k$ involves $R_k$, $G_k$ and $G_k^\prime$ and
\begin{equation}
\label{1.6}
P(\eta) = \int_{\eta_i}^\eta \frac{1}{G_0 (\eta')} d\eta'~,~ \text{for some~} 
\eta_i \in \mathcal{D}
\end{equation}
while $\omega_1 = \e^{i 2 \pi/3}$, $\omega_2 = \e^{-i 2 \pi/3}$ and $\omega_3 =
1$ (the three cubic roots of unity).  In (\ref{1.5}) the limits of
integration satisfy $\eta_{j,k} \equiv \beta_k^{4/9} \chi_j$, where $\chi_1$,
$\chi_2$ and $\chi_3$ are as defined in the last subsection.  The choice of
the domain $\mathcal{D}_k $ for $\eta$ is subject to the conditions below.
\begin{enumerate}
\item{} $\mathcal{D}_k $ must contain a region $\mathcal{S}_0$ that winds around
  $\eta_s$, near $\arg \eta = - \frac{4}{9} \pi$, 
  excluding a $O(\eta_s^{-5/4})$ neighborhood of $\eta_s$ itself. 
  Since the singularity ${\hat \eta}_s$ of $G_0$ is within 
$\e^{-C |\eta_s|}$
  distance of $\eta_s$ and a singularity of $G_0$ is the only point where
  $G_0 =0$, this condition ensures a lower bound on $G_0$ and provides the contour integration $\oint_{S_{\epsilon_1}} $ in
  Remark 2.

\item{} Any point $\eta \in \mathcal{D}_k $ can be connected to
  $\eta_{j,k}$ along a contour $\mathcal{C}_j$ that lies entirely
  within $\mathcal{D}_k$ such that $\Re\left [ \omega_j P
  \right ]$ is increasing monotonically when  the points $\eta_{j,k}$ are
  approached. The integration contour $\mathcal{C}_j$ in (\ref{1.5})
 is chosen to be such a path. Monotonicity
  ensures there is no exponential growth in $k$ ($\beta_k$).
  A fortiori, the points $\eta_{j,k}$, as defined above, are points of
  maximum of $\Re\left [ \omega_j P \right ]$ in
  $\mathcal{D}_k$.  \item{} We must have for $k \ge k_0$, the property
  $\mathcal{D}_{k+1} \cup \mathcal{E}_{k+1} \subset \mathcal{D}_k \cup
  \mathcal{E}_k $. For $k \le k_0$, all $G_k$ can be determined
  through the representation (\ref{1.4.13.0}) on a common domain
  $\mathcal{D}_{k_0} \cup \mathcal{E}_{k_0}$.  The necessity of this
  condition comes from the fact that the $R_k$, needed to determine
  $G_k$ in the domain $\mathcal{D}_k \cup \mathcal{E}_k$, involve
  $G_1$, $G_2$, ..., $G_{k-1}$.
  
\item{} For any $k$, $\mathcal{D}_k \cup \mathcal{E}_k $ must contain the
  region $\mathcal{S}_0$ around the singularity $\eta_s$.  Also,
  for large enough $|\eta|$ in this domain, we must have $\arg \eta \in \left
  [ -\frac{2}{9} \pi + \delta, \frac{2}{9} \pi - \delta \right ]$.  We shall
  furthermore  ensure a nonempty common domain $\mathcal{D} =
  \cap_{k\ge k_0} [\mathcal{D}_k \cup \mathcal{E}_k] $.
  
\item{} To simplify the analysis, we seek domains so that $\mathcal{D}_k \cup
  \mathcal{E}_k$ does not contain turning points (occurring when $\arg \eta =
  \pm\frac{2}{9} \pi $) of  the WKB solutions for large $k$ in
  (\ref{1.4.13}),  $\e^{\beta_k P_j (\chi(\eta)) + W_j (\chi (\eta))}$ (see
  \S4).

\end{enumerate}

\section{Properties of $P(\eta)$ and choice of the  domains $\mathcal{D}$ and $\mathcal{D}_k$}

This section is devoted to the construction of the domains
$\mathcal{D}$ and $\mathcal{D}_k$ corresponding to a particular
$\eta_s$, determined from (\ref{0.3.2.1}) for large enough ${\hat
  n}\in \NN$.  The domains clearly depend on the choice of $\eta_s$.
The construction is relatively involved since monotonicity of
$\Re\left [\omega_j P ~\right ]$ must be ensured, while $P(\eta)$ is
only implicitly known through (\ref{0.3.1})--(\ref{0.3.2}) and
(\ref{1.6}).  Also, such a domain has to contain an annular region
around singularity $\eta_s$.

\begin{Remark}
  In this section, constants such as ${C}$,
    ${K}$, ${\delta}$, $r_i$, $r_0$, etc are positive and independent of
    $\eta$ and $\eta_s$.
\end{Remark}
First we define $\mathcal{D}_A$, part of the region where $G_0 (\eta)
\sim ~\eta^{-1/2} $.
\begin{Definition} For small ${\delta} \in \left (0, \frac{\pi}{63} \right )$
we have 
\label{defD_A}
$$
\mathcal{D}_{A_1}= \left \{ \eta: |\eta| > \frac{1}{2} |\eta_s|,
  ~\arg \eta \in \left ( -\frac{2}{9} \pi + \delta , \frac{2}{9} \pi -
    \delta \right ) \right \}
$$
$$ \mathcal{D}_{A_2} = \left \{ \eta: ~\theta = \arg \eta \in 
\left ( -\frac{4}{9} \pi + \delta, -\frac{2}{9} \pi + \delta \right ], 
|\eta_1 (\theta))|  > |\eta| > \frac{1}{2} |\eta_s|, \right \}
$$
where $\eta_1 (\theta (r)) = 
2 |\eta_s| \e^{-i \frac{4}{9} \pi + i \delta } + r \e^{-i \frac{\pi}{9}} $ for
$  r \ge 0$. We also define $M = |\eta_1 \left ( -\frac{2}{9} \pi + \delta \right )|$ and 
$$ \mathcal{D}_{A} = \mathcal{D}_{A_1} \cup \mathcal{D}_{A_2} $$
(See Figure 1).

\end{Definition}

\begin{Lemma}
\label{D_Adescent}
For any point $\eta \in \mathcal{D}_A $, there exist three piecewise
smooth paths from $\eta$
to $\infty$, $\tilde\eta:={\tilde \eta}_j$, for $j=1,2,3$ , contained in $\mathcal{D}_{A}$ so that on any smooth
segment we have
$$ \frac{d}{ds} \Re \left (\omega_j P [{\tilde \eta} (s) ]\right ) > 
C |{\tilde \eta} (s)|^{1/2} > 0 $$
where $s$ is the arclength. Furthermore, 
$$ | {\tilde \eta} (s) | > C_1 |\eta| > 0 $$
\end{Lemma}

\begin{proof} On the line segment
  ${\tilde \eta} (s) = \eta_0 + s \e^{i \phi} $ in $\mathcal{D}_A$,
   (\ref{1.6}) and largeness of $|{\tilde \eta}(s)|$ and
    $|\eta_s|$, together with the asymptotic behavior $G_0 ({\tilde
      \eta}) \sim {\tilde \eta}^{-1/2}$, imply
$$ \frac{d}{ds} \ \Re( \omega_j P) = 
\ \Re\left [ \frac{\omega_j}{G_0 ({\tilde \eta} (s))} {\tilde \eta}^\prime (s) \right ] 
~\sim ~|{\tilde \eta}(s)|^{1/2} \cos \left (\theta_j + \frac{\theta_{\tilde \eta}}{2} + \phi \right )
$$
where $\theta_j = ~\arg \omega_j \in \left \{ \pm{\frac{2}{3}}
  \pi , 0 \right \}$, $\theta_{\tilde \eta} = \arg {\tilde \eta}(s)$.
For suitable  $\phi$ and $\eta_0$, it is easy to see that for
any $\eta \in \mathcal{D}_A$ and $j=1,2,3$, there exists a polygonal line
so that $\cos \left (\theta_j + \frac{1}{2} \theta_{\tilde \eta} + \phi
\right )> C > 0$.  Further, the line can be chosen so that
$|{\tilde \eta } (s)| \ge |\eta|$.
\end{proof}

\begin{Definition}
\label{defL}
$$ L = \left \{ \eta : \arg \eta = -\frac{4}{9} \pi + \delta, 
\frac{1}{2} |{\eta_s}| < |\eta| < 2 |\eta_s | \right \} $$
\end{Definition}

\begin{Definition}
Let $\mathcal{S}_0$ be a region around $\eta_s$ 
(the singularity of $U(\zeta(\eta))$ 
in (\ref{0.3.2})) defined by  
\begin{equation}
\label{2.2}
\mathcal{S}_0 = \left \{ \eta:  ~r_i |\eta_s|^{-5/4} < |\eta-\eta_s| <  r_o |\eta_s|^{-5/4}, 
~\arg (\eta-\eta_s) 
~\in ~\left ( -\pi +\frac{\pi}{18}, \frac{\pi}{18} + \pi \right )  \right \} 
\end{equation}
with $0 < r_i < r_0$,  small enough to ensure 
$\arg U_0^{1/2} \in \left [-\frac{2}{5} \pi, \frac{2}{5} \pi \right ]$ (see relation 
(\ref{2.2.6.8})).
\end{Definition}

\begin{Definition}
  We define 
\label{DefDT}
\begin{multline*}
\mathcal{D}_{T,1} = \left \{ \eta: |\eta^{3/2} - \eta_0^{3/2} | < 3
\left (\frac{1+{B}_0}{1-{B}_0}\right )^2 r ~~{\rm for} ~
\eta_0 \in \mathcal{S}_0  ~{\rm and} \right .\\
\left . \left |\frac{1-\sqrt{U}}{1+\sqrt{U}} \right | < {B}_0
\e^{-K_4 |\eta_s|^{3/4}r} < 1, ~{\rm for~some}~r \in 
[0, |\eta_s|^{1/4}] 
\right \}
\end{multline*}
For $r \in \left [ 0, \frac{2 \delta}{K_3} |\eta_s|^{3/2} \right ]$,
we define
\begin{align*}
\mathcal{D}_{T2, r}:= \Bigg \{ \eta\,{\text{\Large \rm :}} \ \ \ &(1)\  \arg \eta 
\in\left( \arg \eta_s - K_5|\eta_s|^{-5/4} +
      K_3 |\eta_s|^{-3/2} r,-\frac{4}{9} \pi +
      \delta,\right] \\
    &(2)\   |\eta|^{3/2} \in\left(|\eta_s|^{3/2} - K_5 |\eta_s|^{1/4} - 3r ,
   |\eta_s|^{3/2} + K_5 |\eta_s|^{1/4} + 3 r  \right) \\
   &(3)\   |U-1| < 
5 {B}_0 \e^{-K_4 |\eta_s|} \e^{-K_1 |\eta_s|^{3/4} r} \Bigg\} 
 \end{align*}
\begin{align*}
\mathcal{D}_{T2}:=& \bigcup_{r \in I} \mathcal{D}_{T2,r}, ~~{\rm where} ~ 
~I = \left [ 0, \frac{2 \delta}{K_3} |\eta_s|^{3/2} \right ]\\
\mathcal{D}_T:=& \mathcal{D}_{T,1} \cup \mathcal{D}_{T,2}
~~~~~~~~~~~~~~~~~~~~~~~~~~~~~~~~~~~~~~~~~~~~~~~~~~~~~~~~~~~~~~~~~~~~~~~~~~\\\
\mathcal{D}:=&\mathcal{D}_A\cup\mathcal{D}_T
\end{align*}
\end{Definition}
\begin{Theorem}
\label{ThmDT}
For large $|\eta_s|$, for any point $\eta \in \mathcal{D} $, there
exist ${B}_0$, $K_i$ and piecewise smooth paths from $\eta$ to
$\infty$ $\tilde\eta_j=:\tilde{\eta}$ ($j=1,2,3$) contained in
$\mathcal{D}$, so that on any smooth subsegment we have
\begin{equation}
\label{monoton}
\frac{d}{ds} \Re\left (\omega_j P [{\tilde \eta} (s) ]\right )
 > C |{\tilde \eta} (s)|^{1/2} > 0 
\end{equation}
Furthermore
\begin{equation}
\label{ratio}
| {\tilde \eta} (s) | > C_1 |\eta| > 0 
\end{equation}
\end{Theorem}
  For the proof, given at the end of \S2, we need a few more
  definitions, constructions and lemmas.
\begin{Remark} 
\label{remThmDT}
Note that by Lemma \ref{D_Adescent}, it is enough to show that
for any $\eta \in \mathcal{D}_T $, we can choose a path for each of
$j=1,2,3$ connecting $\eta$ to $\eta_L \in L$, entirely within
$\mathcal{D}_{T} $ so that the monotonicity property (\ref{monoton})
is satisfied. Noting also that since the ratio of any two values
of $\eta \in \mathcal{D}_T$ is bounded by a constant independent of
$\eta_s$, the second part of Theorem \ref{ThmDT} follows.
\end{Remark}

\begin{Definition}\label{defDk}
  For $k \ge k_0$, we define
  $$\eta_{1,k}= \eta_{3,k} + \rho_0 \e^{-i \frac{2}{3} \pi} ,
  \eta_{2,k} = \eta_{3,k} + \rho_0 \e^{i \frac{2}{3} \pi } ,
  \eta_{3,k} = \epsilon \beta_k^{4/9} $$
  where $\rho_0$ is chosen
  so that $\arg \eta_{1,k} =-\frac{2 \pi}{9} -\delta $, $\arg
  \eta_{2,k} =\frac{2 \pi}{9} +\delta $ for $0
  <\delta <\frac{\pi}{63}$. The parameter $\epsilon$ is small, but
    independent of $k$, as needed in Lemma \ref{L3.5}, and $k_0$
  is chosen large enough so that for $k \ge k_0$, we have  $\epsilon
  \beta_k^{4/9} > M $, for $M$ as defined in Definition \ref{defD_A}. 
  We define a boundary ${\partial E}_k =
  {\partial E}_k^- \cup {\partial E}_k^+ $ where ${\partial E}_k^-$ is
  the straight line joining $\eta_{3,k}$ with $\eta_{1,k}$ and
  ${\partial E}_k^+$ is the straight line joining $\eta_{3,k}$ to
  $\eta_{2,k}$.  We then define $\mathcal{D}_k$ (See Fig. 2)
$$\mathcal{D}_k 
= \mathcal D \backslash \mathcal{E}_k 
$$ 
\end{Definition}

\begin{Lemma}
\label{monEk}
Given $j =1,2$ or $3$, for any $\eta \in {\partial E}_k $, 
the path ${\tilde \eta} (s)$ from $\eta$ to $\eta_{j,k}$ along
${\partial E}_k$ satisfies the monotonicity property (\ref{monoton}). 
\end{Lemma} 

\begin{proof}
On ${\partial E}_k^+$ we note that 
$$  \frac{d}{d\rho} \Re\left [ \omega_3 P (\eta_{3,k} 
+ \rho \e^{i 2 \pi/3 } ) \right ]
\sim\Re\left [ \e^{i 2 \pi/3} \eta^{1/2} \right ]
=-|\eta|^{1/2} \sin \left (\frac{\pi}{6} +\frac{1}{2} \arg\eta \right )<- C |\eta|^{1/2}$$
for some positive constant $C$. 
On ${\partial E}_k^-$ we note that 
$$
\frac{d}{d\rho} \Re\left [ \omega_3 P (\eta_{3,k} +
  \rho \e^{-i 2 \pi/3 } ) \right ] \sim \Re\left [ \e^{-i 2
    \pi/3 } \eta^{1/2} \right ] =-|\eta|^{1/2} \sin \left
  (\frac{\pi}{6} +\frac{1}{2} \arg\eta \right ) <-C
|\eta|^{1/2}$$
for some positive constant $C$. It is therefore
clear that the path ${\tilde \eta}_3 (s)$ from ${\tilde \eta}$ to
$\eta_{3,k}$ satisfies the monotonicity property (\ref{monoton}). In a
similar manner, it is seen that $\Re \left [ \omega_1
  P \right]$ and $\Re\left [ \omega_2 P \right ]$ satisfy
(\ref{monoton}) on a path from $\eta $ to $\eta_{j,k}$ for $j=1,2$
along ${\partial E}_k$.
\end{proof}

The lemma above, together with Theorem \ref{ThmDT} prove the following Corollary:

\begin{Corollary}
\label{Cor:Property1}
({\bf Property 1:}) For all sufficiently large $k$, given any point $\eta \in
\mathcal{D}_k$, there exists a  piecewise smooth path  $\mathcal{C}_j$ for each $j=1,2,3$
from $\eta$ to $\eta_{j,k} $ such that the path is entirely in
$\mathcal{D}_k$ and
\begin{equation}
  \label{Property1}
  \frac{d}{ds} \Re\left (\omega_{j} P ({\tilde \eta} (s) 
\right ) \ge C
|\tilde \eta|^{1/2}>0
\end{equation}
Furthermore, if ${\tilde \eta} \in
\mathcal{C}_j$,   we then have $|\tilde \eta| > C |\eta| $.  
\end{Corollary}

\begin{Remark} To prove Theorem \ref{ThmDT}, we introduce
   three autonomous flows as follows.
\end{Remark}

\begin{Definition}
  Let $\eta_j (t, \eta_0)$ be the solution to the differential 
equation
\begin{equation}
\label{2.1}
\dot{\eta} = i \e^{-i \phi_j} \omega_j^{-1} G_0(\eta),~{\rm where} ~\omega_1 = \e^{i 2 \pi/3}, 
~\omega_2 = \e^{-i 2 \pi/3},~\omega_3 = 1   
\end{equation}
with initial condition $\eta_j (0, \eta_0) = \eta_0$, where
$\phi_j$ are given by:
\begin{equation}
\label{phiset1}
\phi_1 = \frac{\pi}{3}, ~\phi_2 = \frac{6 \pi}{7} , ~\phi_3 = \frac{2}{3} \pi 
\end{equation}
\end{Definition}
\begin{Remark}
\label{rem9}
We note from (\ref{1.6}) that for any choice $\phi_j \in (0, \pi) $, 
$$ \frac{d}{dt} ~\Re\left [ \omega_j P (\eta_j (t, \eta_0)) \right ]
= ~\cos \left (\phi_j - \frac{\pi}{2} \right ) > 0 $$
Hence using arclength parameterization we have
$$ \frac{d}{ds} ~\Re\left [ \omega_j P ({\tilde \eta} (s) \right ] = 
\frac{\cos \left (\phi_j - \frac{\pi}{2} \right )}{|G_0 ({\tilde \eta(s)}) |} 
>C |{\tilde \eta} (s) |^{1/2} $$
when ${\tilde \eta} (s) \in \mathcal{D}$.
Thus, the differential equation (\ref{2.1}) generates 
ascent paths for $\Re [\omega _j P]$. 
\end{Remark} 

\begin{Lemma} 
\label{S0subsetD}
There exists a ${B}_0$ so that 
$\mathcal{S}_0 \subset \mathcal{D}_T $. 
\end{Lemma} 
\begin{proof}
Since for $\eta \in \mathcal{S}_0$, the corresponding $U$ determined from
(\ref{0.3.2}) has upper and lower bounds independent of $\eta_s$, as discussed
in Remark \ref{U01}. Also, from (\ref{2.2.6.8}), for $\eta \in \mathcal{S}_0$,
$\arg U^{1/2} \in \left [-\frac{2}{5} \pi, \frac{2}{5} \pi \right ]$.  Thus,
it follows that for $\eta \in \mathcal{S}_0$,  we have  $
|1-\sqrt{U}|/|1+\sqrt{U}| < {B}_0 $ for some ${B}_0 < 1$.  Thus, for
some  ${B}_0$, we have $\mathcal{S}_0 \subset \mathcal{D}_{T,1} \subset
\mathcal{D}_{T}$.
\end{proof}

\begin{Definition}
It is convenient to define,  see (\ref{2.1}) and (\ref{phiset1}),
$$\nu_j =i \e^{-i \phi_j} \omega_j^{-1} $$ 
\end{Definition}

\begin{Remark}
 It follows that
\begin{equation} 
\label{2.2.9}
\arg \nu_1 =-\frac{\pi}{2},~~\arg \nu_2 =\frac{13 \pi}{42},~~\arg \nu_3
 =-\frac{\pi}{6} ~~
\end{equation}
The specific choice of $\phi_j$ (and thus of $\nu_j$) is unimportant, but it
is essential that $\phi_j $, $\arg \nu_j$ remain in compact
subintervals of $(0, \pi)$ and $\left ( -\frac{2}{3} \pi, \frac{\pi}{3} \right
)$ respectively, independent of $\eta_s$ and ${\delta}$.

\end{Remark}

In order to study the solution to (\ref{2.1}) 
near $\eta_s$, 
it is convenient to think of $U(t) = U(\eta(t))$ as an
unknown together with $\eta(t)$. 
Using (\ref{0.3.1.5}), (\ref{0.3.2}) and (\ref{2.1}), it follows
that 
\begin{equation}
\label{2.2.6.3}
\frac{2}{3} \frac{d}{dt} \eta^{3/2} = \nu_j U \left [ 1 + E_1 (\eta) \right ]
~~{\rm where} ~~E_1 (\eta) =\frac{\eta^{1/2} G_0 (\eta) - U}{U} 
\end{equation}
\begin{multline}
\label{2.2.6.4}
\frac{d}{dt} U = 
-\alpha_j |\eta_s|^{3/4} \sqrt{U} (U-1)  
\left [ 1 + E_1  \right ] \left [ 1 + E_2 \right ], 
~{\rm where} ~\\
E_2 (\eta) =\left [ \frac{\nu_j}{\alpha_j |\eta_s|^{3/4}} \left ( 
\frac{i \sqrt{2}}{3} \eta^{3/4} + \frac{9}{8 \eta^{3/2}} \right ) - 1 \right ]
\end{multline}
where $E_1 $, $E_2$ will be shown to be small for large $|\eta_s|$
in the range of integration  and
$$ \alpha_j = 
\frac{\nu_j}{|\eta_s|^{3/4}} \left [ \frac{i \sqrt{2}}{3} \eta_s^{3/4} + \frac{9}{8 \eta_s^{3/2}} \right ] $$
The initial condition $U_0$ satisfies
\begin{equation}
\label{2.2.6.6}
\frac{i 4\sqrt{2}}{27} \eta_0^{9/4}-\frac{i 4\sqrt{2}}{27} \eta_s^{9/4}
 + \frac{9}{8} \ln \left (\frac{\eta_0}{\eta_s}\right) 
=-\ln \frac{1 - \sqrt{U_0}}{1+\sqrt{U_0}} - 2 \sqrt{U_0}
\end{equation}

\begin{Remark} 
It is to be noted that with $\phi_j$ given by (\ref{phiset1}) and using the
fact that as ${\hat n} \rightarrow \infty$ (i.e. as $|\eta_s| \rightarrow
\infty$), we get $\arg\eta_s \rightarrow -\frac{4 \pi}{9} $. It follows that in this
limit,
\begin{equation} 
\label{2.2.7}
\arg\alpha_1 \rightarrow -\frac{\pi}{3},~
\arg\alpha_2 \rightarrow \frac{10 \pi}{21} ,~\arg\alpha_3
 \rightarrow 0 
\end{equation}
It is important for us that  $\arg \alpha_j\in \left (-\frac{\pi}{2}, \frac{\pi}{2}\right ) $.
\end{Remark}

\begin{Lemma}
\label{L0.2.2.0}
For suitable $ K_i$  and $\delta < \frac{K_3}{24}$, 
if $\eta_{0,0} \in \mathcal{D}_{T_2}$, then for sufficiently large $|\eta_s|$
and
some $t\in\left (0, \frac{2 \delta }{K_3}
  |\eta_s|^{3/2} \right )$, $\eta_j (t; \eta_{0,0} )$, leaves
$\mathcal{D}_{T_2}$ through $L$.
\end{Lemma}

\begin{proof} 
The differential equation  
satisfied by $\eta$ and the corresponding $U$ for 
$0 \le t \le 2 \frac{\delta}{K_3} |\eta_s|^{3/2} $ 
is given by:
$$\frac{2}{3}\frac{d}{dt}  \eta^{3/2} = \nu_{j} [ 1 + E_3 ] ~~;~~\frac{d}{dt} U = - \nu_{j}
\left ( \frac{i \sqrt{2}}{3} \eta^{3/4} + \frac{9}{8\eta^{3/2}} \right ) (U-1)
                    (1+ E_4 ) $$
where 
$$ E_3 = E_1 U + (U-1) ,~~E_4 = (\sqrt{U} - 1) + \sqrt{U} E_1 ,  $$
where $E_1$, $E_2$ are  defined in (\ref{2.2.6.3}) and (\ref{2.2.6.4}). It follows that
\begin{equation}
\label{eq32}
\eta^{3/2} = \eta_{0,0}^{3/2} 
+ \frac{3}{2} \nu_{j} \int_0^t (1+E_3 )dt 
\end{equation}
\begin{equation}
\label{eqUm1}
(U-1) = (U_{0,0} - 1) 
\exp \left [ - \nu_{j} \int_0^t \left ( \frac{i \sqrt{2}}{3} \eta^{3/4} + 
\frac{9}{8\eta^{3/2}} \right ) ( 1 + E_4 ) dt \right ] 
\end{equation}
where $U_{0,0}$ is obtained from (\ref{2.2.6.6}) by substituting $\eta_0 = \eta_{0,0}$.
It is convenient to define the  leading order equations
\begin{equation}
\label{eq32.tilde}
{\tilde \eta}^{3/2} = \eta_{0,0}^{3/2} 
+ \frac{3}{2} \nu_{j} t ~~;~~
{\tilde U}-1 = (U_{0,0} - 1)
\exp \left [ - \nu_{j} \int_0^t \left ( \frac{i \sqrt{2}}{3} {\tilde \eta}^{3/4} + 
\frac{9}{8{\tilde \eta}^{3/2}} \right ) dt \right ] 
\end{equation}
From (\ref{eq32}), (\ref{eqUm1}) and (\ref{eq32.tilde})
it follows that
\begin{equation}
\label{eq32net}
\eta^{3/2} - {\tilde \eta}^{3/2} = \frac{3}{2} \nu_{j} \int_0^t E_3 dt 
\end{equation}
\begin{multline}
\label{eqUm1net}
U-{\tilde U} = (U_{0,0}-1) \Bigg \{ 
\exp \left [ - \nu_{j} \int_0^t \left ( \frac{i \sqrt{2}}{3} \eta^{3/4} + 
\frac{9}{8\eta^{3/2}} \right ) ( 1 + E_4 ) dt \right ] 
\\- \exp \left [ - \nu_{j} \int_0^t \left ( \frac{i \sqrt{2}}{3} {\tilde \eta}^{3/4} + 
\frac{9}{8{\tilde \eta}^{3/2}} \right ) dt \right ] \Bigg \}
\end{multline}
From (\ref{eq32.tilde})  it follows that
$$ \arg \eta_{0,0} + 3 t |\eta_s|^{-3/2} \ge \arg {\tilde \eta} \ge \arg \eta_{0,0} 
+ \frac{3}{2} K_3 t |\eta_s|^{-3/2}
~~~{\rm where} ~2 K_3 = \min_{j=1,2,3} \sin \left ( \frac{2}{3} \pi + 
\nu_j \right )
$$ 
$$|{\tilde \eta}|^{3/2} \in 
\left ( |\eta_{0,0}|^{3/2} - \frac{3}{2} t, 
|\eta_{0,0}|^{3/2} + \frac{3}{2} t \right )$$
Using these relations in (\ref{eq32.tilde}) we have
$$ |({\tilde U} - 1)| = |(U_{0,0} - 1)|  \e^{-2 K_1 |\eta_s|^{3/4} t } $$
where $2 K_1$ is a lower bound (independent of $\delta$) of 
$$ \Re \left \{ \frac{\nu_j}{2 |\eta_s|^{3/4}}  \left [ \frac{i \sqrt{2}}{3} {\tilde \eta}^{3/4} + 
\frac{9}{8 {\tilde \eta}^{3/2}} \right ] \right \} $$
for $\tilde \eta$ restricted to the  domain $|\tilde \eta| > \frac{1}{2} |\eta_s|$, 
$\arg \tilde \eta \in \left [ -\frac{4}{9} \pi - \delta, -\frac{4}{9} \pi + 4 \delta \right ]$. 
Thus, for some $t$ in 
$ 0 \le t \le  \frac{3 \delta}{2 K_3} |\eta_s|^{3/2} $, ${\tilde \eta}$
leaves
the domain $\mathcal{D}_{T_2}$ through the segment of $L$, when
$\frac{15}{32} |\eta_s | < |{\tilde \eta}| < \frac{3}{2} |\eta_s| $. 

Now, we show that $\eta$ is close to
${\tilde \eta}$ and hence has roughly the same behavior. 
We define
$$ (\zeta, V) = \left ( \frac{\eta^{3/2} - {\tilde \eta}^{3/2}}{t_1}, 
\frac{U - {\tilde U}}{U_{0,0}-1} \right )$$
on the interval $[0, t_1]$, for $0 < t_1 \le \frac{2 \delta}{K_3} |\eta_s|^{3/2} $.
We introduce the norm 
$$ \| (\zeta, V) \|_\infty = \sup_{ 0 \le t \le t_1} |\zeta (t) |
+ \sup_{ 0 \le t \le t_1} \e^{\frac{3}{2} K_1 |\eta_s|^{3/4} t} |V (t) | $$
and consider the right side of (\ref{eq32net}) and (\ref{eqUm1net}) as
the mapping 
$$\left ( \mathcal{F}_1 (\zeta, V), \mathcal{F}_2 (\zeta, V) \right ) $$
of the ball
$$
\mathcal{B} = 
\left \{ (\zeta (t), V (t) ) : \| (\zeta, V) \|_\infty < \epsilon_1 \right \}
$$ 
for  some small $\epsilon_1$ 
in the Banach space of pair of continuous functions
$(\zeta (t), V(t) ) $ of $t$ in the
interval $[0, t_1]$ for $t_1 < \frac{2 \delta}{ K_3} |\eta_s|^{3/2}$.  

Using  the  smallness of $E_3$ and $E_4$ for large $\eta$ 
it can be checked directly that
$$ \left ( \mathcal{F}_1 (\zeta, V), \mathcal{F}_2 (\zeta, V) \right )
\in \mathcal{B} $$ 
and that 
$$ \| \left ( \mathcal{F}_1 (\zeta_1, V_1 ) , \mathcal{F}_2 (\zeta_1, V_1) 
\right )
 -\left ( \mathcal{F}_1 (\zeta_2, V_2 ) , \mathcal{F}_2 (\zeta_2, V_2) 
\right ) \|_\infty \le \epsilon_2 \| (\zeta_1, V_1) - (\zeta_2, V_2 ) \|_\infty $$
for some $\epsilon_2 < 1$ and the map is contractive.
Thus, there is a unique solution to the integral 
system (\ref{eq32net}) -- (\ref{eqUm1net})
for $(\zeta (t) , V (t))$ in $\mathcal{B}$. 
In particular, this implies that
\begin{equation}
\label{resUeta}
| U(t) - 1 | \le |U_{0,0} - 1| \e^{- K_1 |\eta_s|^{3/4} t} ,~
   | (\eta (t))^{3/2} - \eta_{0,0}^{3/2} | \le 3 t   
\end{equation}
Hence, with $r$ as in  the
definition of $\mathcal{D}_{T_2, r}$ we have
$$ |U - 1| \le 5 B_0 \e^{-K_4 |\eta_s| } 
\e^{-K_1 |\eta_s|^{3/4} (t +r )} $$
$$ \arg \eta \ge \arg \eta_{0,0} + K_3 |\eta_s|^{-3/2} t
\ge \arg \eta_s -K_5 |\eta_s|^{-5/4} + K_3 |\eta_s|^{-3/2} (t+r) $$
$$ |\eta^{3/2}| \in \left ( |\eta_s|^{3/2} -K_5 |\eta_s|^{1/4} - 3 (t +r), 
|\eta_s|^{3/2} + K_5 |\eta_s|^{1/4} + 3 (t +r) \right ) $$ 
Therefore, from the definition of $\mathcal{D}_{T,2}$, for
small enough $t+r$,  we have $\eta \in \mathcal{D}_{T,2}$, while from continuity, there
exists some larger 
$t+r  \le \frac{2 \delta}{K_3}|\eta_s|^{3/2} $ 
for which $\eta \in L$ as it exits $\mathcal{D}_{T_2} $. 
\end{proof}

\begin{Lemma}
\label{L0.2.1.3}
Let $\eta_{0,0} \in \mathcal{D}_{T,1}$.
Define
$$
{\hat \eta} = \eta_{j} (t,\eta_{0,0}) $$
Then, there exist ${B}_0$ and
${K}_i$ so that ${\hat \eta} \in \mathcal{D}_{T_1} \cup
\mathcal{D}_{T_2}$ for large $|\eta_s|$ and $ 0 \le t \le
|\eta_s|^{1/4}$.
\end{Lemma}

\begin{proof} Note that for any 
$0 \le t \le |\eta_s|^{1/4} $, we write (\ref{2.2.6.3})
and (\ref{2.2.6.4}) as
\begin{equation}
\label{2.2.8.1}
\eta^{3/2} = \eta_{0,0}^{3/2} + \frac{3}{2} \int_0^t \nu_{j} U (1+ E_1) dt ~,~
\frac{1-\sqrt{U}}{1+\sqrt{U}} = b_{0,0} 
\exp \left \{-\alpha_{j} |\eta_s|^{3/4} 
\int_0^t (1+E_2) (1+ E_1 )dt \right \} 
\end{equation}
where $|b_{0,0}| < B_0 \e^{-K_4 |\eta_s|^{3/4} r}$, with
$B_0$ chosen in accordance to Lemma 
\ref{S0subsetD} and 
$ 2 K_4 := \min_{j} \cos \left ( \frac{\pi}{6} + \nu_j \right ) $.
We introduce ${\tilde \eta} (t)$ and
${\tilde U} (t)$ (describing leading behavior) by
$$ \frac{1-\sqrt{\tilde U}}{1+\sqrt{\tilde U}} 
= b_{0,0} \e^{-\alpha_{j} |\eta_s|^{3/4} t} ~~,~~
{\tilde \eta}^{3/2} = \eta_{0,0}^{3/2} +
\frac{3}{2} \nu_{j} \int_0^t {\tilde U} (t') dt' $$
It is to be noted that
\begin{equation}
\label{etatilde1}
\eta^{3/2} - {\tilde \eta}^{3/2} = \frac{3}{2} \int_0^t \nu_j [(U-1) + U E_1 ] dt 
\end{equation}
\begin{equation}
\label{etatilde2}
\frac{1-\sqrt{U}}{1+\sqrt{U}} - \frac{1-\sqrt{\tilde U}}{1+\sqrt{\tilde U}} 
= b_{0,0} 
\e^{
-\alpha_j t |\eta_s|^{3/4}} \left [ 
\exp \left \{-\alpha_{j} |\eta_s|^{3/4} 
\int_0^t [(1+E_2) (1+ E_1 )-1] dt \right \} - 1 \right ]
\end{equation} 
We note that $\frac{5}{3} K_4$ is a lower bound for $\Re ~[\alpha_j]$ for
$|\eta_s|$ large.   
It is convenient to define the pair of continuous functions, 
$$ (\zeta (t), V(t) ) = 
\left ( \eta^{3/2} (t)-{\tilde \eta}^{3/2} (t),  
\frac{1-\sqrt{U (t)}}{1+\sqrt{U(t)}} - \frac{1-\sqrt{\tilde U(t)}}{1+\sqrt{\tilde U(t)}}  \right )$$
and  the norm
$$ \| (\zeta, V) \|_\infty = \sup_{0 \le t \le t_1} |\zeta(t) | + 
\sup_{0 \le t \le t_1} \e^{\frac{3}{2} K_4 |\eta_s|^{3/4} t} |V(t)|
$$  
for $t_1 \in (0, |\eta_s|^{1/4} )$. 
Consider the right hand side of
(\ref{etatilde1}) and (\ref{etatilde2}) as a mapping 
$\left ( \mathcal{F}_1 (\zeta, V), 
\mathcal{F}_2 (\zeta, V) \right )$ on the ball  
$$ \mathcal{B} = \left \{ (\zeta, V): \|(\zeta, V) \|_\infty < \epsilon_1 t_1 \right \}
$$ 
Using smallness of $E_1$, $E_2$ and their derivatives with respect to
$\eta$, it can be readily checked that 
$\left ( \mathcal{F}_1, \mathcal{F}_2 \right )$ is a contractive mapping of
the ball $\mathcal{B}$ into itself; hence the solution $(\zeta, V)$ satisfying
(\ref{etatilde1}) and (\ref{etatilde2}) is in $\mathcal{B}$
for large $|\eta_s|$.
In particular, since $\Re \alpha_j > \frac{5}{3} 
 K_4 $  we have 
$$ \Bigg\lvert \frac{1-\sqrt{U (t})}{1+\sqrt{U (t)}} 
\Bigg\rvert \le {B}_0  
\e^{- K_4 |\eta_s|^{3/4} (t + r)} ,~
\lvert [\eta(t)]^{3/2} - \eta_{0}^{3/2} \rvert \le 3 
\left ( \frac{1+{B}_0}{1-{B}_0} \right )^2 
(t+r)   $$
There are two cases: if $t + r \le |\eta_s|^{1/4}$, then clearly $\eta \in \mathcal{D}_{T_1}$.
If $|\eta_s|^{1/4} \le t + r \le 2 |\eta_s|^{1/4} $, from the definition of $\mathcal{D}_{T_2}$, 
it follows $\eta \in \mathcal{D}_{T,2}$, 
with $K_5 = 6 \left ( \frac{1+B_0}{1-B_0} \right )^2$.
\end{proof}

\noindent{\bf Proof of Theorem \ref{ThmDT}}. 
From Lemmas \ref{D_Adescent},
\ref{L0.2.1.3}, \ref{L0.2.2.0} (see Remark \ref{remThmDT} as well),
it is clear that  
that the domain
$\mathcal{D} = \mathcal{D}_{T} \cup \mathcal{D}_A$ is invariant
under the flows ${\tilde \eta}_j (s)$. From Remark \ref{rem9},
Theorem \ref{ThmDT} follows.    

\section{Properties of $P_j (\chi)$ and choice of domain $\mathcal{\mathcal{E}}$}

\begin{Remark} 
  The WKB solution for large $\beta_k$ of the homogenous  equation
  $\mathcal{L}_k u = 0$ (see \S\ref{S:WKB}, item 5) is not uniformly valid in the domain $\mathcal{D}$ for
  large $\eta$.  To invert the operator $\mathcal{L}_k$ in the   regime
  $\eta= O(\beta_k^{4/9})$, we introduce the scaled variables:
\begin{equation}
\label{3.1}
\chi = \beta_k^{-4/9} \eta 
\end{equation}
The WKB solution to the homogeneous equation is then of the form
\begin{equation}
\label{3.2}
\e^{\beta_k P_j (\chi) + W_j (\chi)}  
~{\rm where} ~\alpha = P_j^\prime ~~{\rm are ~roots ~of~the~cubic} ~\alpha^3  
+ \frac{2}{9}\alpha \chi^{5/2} - \chi^{3/2} = 0 
\end{equation}
We now choose a domain $\mathcal{E}$ where the WKB solution is valid.
First, we define a boundary ${\partial \mathcal{E}}_L$, 
which corresponds in the $\chi$ plane to ${\partial E}_k$ (see Definition
\ref{defDk}).
\end{Remark}

\begin{Definition}
  Let ${\partial \mathcal{E}}_L = \left \{ \chi: \eta = \beta_k^{4/9} \chi \in
    {\partial E}_k \right \} $.  We define ${\partial \mathcal{E}}_L^+$ and
  ${\partial \mathcal{E}}_L^-$ analogously in terms of ${\partial E}_k^+$ and
  ${\partial E}_k^-$, (see Definition \ref{defDk}).
\end{Definition}

\begin{Definition}
\label{D3.2}  We let
$$
\mathcal{E} ~= ~ \left \{ \chi : \chi ~{\rm to~the~right~of}~
  {\partial E}_L, ~~\arg\eta \in ~\left [ -\frac{2 \pi}{9} + \delta ,
    ~\frac{2\pi}{9} -\delta \right ] \right \} $$
(See Fig. 3.)
It is also convenient to
define
$$ \mathcal{E}_k = 
~\left \{ \eta : \beta_k^{-4/9} \eta = \chi ~\in ~\mathcal{E} \right \} $$
\end{Definition}

\begin{Remark}
  Note that for  large $k$  we have the following properties :
  $\mathcal{D} \subset \mathcal{D}_k \cup \mathcal{E}_k $ and 
  $\mathcal{D}_{k+1} \cup \mathcal{E}_{k+1} \subset \mathcal{D}_k \cup
  \mathcal{E}_k $.  This follows from the 
  construction of $\mathcal{D}_k$
  and $\mathcal{E}_k$. Our strategy is to prove the bounds in
  Lemma \ref{L0.2} in the domain $\mathcal{D}_k \cup \mathcal{E}_k$
  based on bounds on all previous $G_j$, $j=1$,2...,$(k-1)$
  established on the domains $\mathcal{D}_j \cup \mathcal{E}_j $  (which contain $\mathcal{D}_k \cup \mathcal{E}_k$). 
    The large $k$
  requirement is not restrictive, since for any fixed $k_0$ it is
  possible to choose $A$ large enough so that the bounds in Lemma
  \ref{L0.2} hold for $ 1 \le j \le k_0$.
  
\end{Remark}    

\noindent The main theorem in this section is the following.

\begin{Theorem} 
\label{T3.1}
For any $\chi \in \mathcal{E}$, it is possible to
choose a path $\mathcal{C}_j $ connecting
$\chi$ to $\infty \e^{i \theta_j}$, 
where $\theta_1 = -\frac{2\pi}{9} + \delta$,
$\theta_2 = \frac{2}{9} \pi -\delta$ 
and $\theta_3 = 0$ so that, except for a finite set of points,
$$ \frac{d}{ds} \Re\left [P_{1,2} ({\tilde \chi} (s) ) \right ] 
~\ge ~ C |{\tilde \chi} (s)|^{5/4} ~>~0 $$
and  
$$ \frac{d}{ds} \Re \left [P_{3} ({\tilde \chi} (s) )\right ] 
~\ge ~ \frac{C}{|{\tilde \chi} (s)|} ~>~0 $$
where $s$ is the arc-length 
increasing towards $\infty$ and
the (different) constants $C$ above are  independent  of $\chi$.
Furthermore, for 
$|\chi|$ sufficiently large in $\mathcal{E}$, 
and with ${\tilde \chi} \in \mathcal{C}_j$
as above, we have $ |{\tilde \chi}| ~>~C ~|\chi|$ for $C ~>~0$ 
independent of $\chi$ and ${\tilde \chi}$.  
\end{Theorem}
\begin{proof} This follows, after a few Lemmas,
at the end of \S3.
\end{proof}

\begin{Remark} Though the domain $\mathcal{E}$ restricts
    the size of $|\chi|$ (it is bounded below), it is convenient to
  first consider the properties of $P_j$ on an enlarged domain
  $\mathcal{E}_0$ with no restriction on $|\chi|$ 
  and larger width:
\end{Remark}

\begin{Definition}
$$\mathcal{E}_0 = ~\left \{ \chi: ~\arg\chi \in \left [ -\frac{2}{9}\pi
,~\frac{2}{9} \pi \right ] \right \} 
$$
\end{Definition}

It is convenient to associate each $P_j$ with a first order
differential equation as follows. Note from (\ref{3.2}) that 
with $P_j^\prime: = \chi^{5/4} \psi$ we have 
\begin{equation}
\label{3.3}
\chi^{-9/4} = \psi^3 + \frac{2}{9} \psi 
\end{equation}
Now, we consider the trajectory in the complex $\chi$ plane generated by the
differential equation
\begin{equation}
\label{3.4}
\frac{d\chi}{dt} = \frac{1}{P_j^\prime (\chi)} ~~{\rm implying}~~  
\frac{4}{9} \frac{d \chi^{9/4}}{dt} = \frac{1}{\psi} 
\end{equation}
The solution with initial value $\chi_0$ will be denoted by $\chi_j (t;
\chi_0)$.  Using (\ref{3.3}), it follows that
\begin{equation}
\label{3.5}
\frac{d\psi}{dt} = -\frac{\psi (2+ 9 \psi^2)^2}{4 (2+ 27 \psi^2)} 
\end{equation}
For large $\chi \in \mathcal{E}$ it is clear from (\ref{3.3}) 
that  the three possible behaviors of $\psi$  are 
$\psi \sim i \sqrt{\frac{2}{9}}$, $\psi = - i \sqrt{\frac{2}{9}} $ and
$\psi \sim ~\frac{9}{2} \chi^{-9/4} $.     
We associate these behaviors with $P_1^\prime$, $P_2^\prime$ and $P_3^\prime$ respectively,
so that 
\begin{equation}
\label{3.5.1}
P_1^\prime ~\sim ~i ~\sqrt{\frac{2}{9}} \chi^{5/4} ~,~     
P_2^\prime ~\sim ~-i ~\sqrt{\frac{2}{9}} \chi^{5/4} ~,~ 
P_3^\prime ~\sim ~\frac{9}{2} \chi^{-1} 
\end{equation}

\begin{Remark} 
\label{R3.1}
Note that $ds = \lvert \frac{d{\tilde \chi}}{dt} \rvert dt
$, and so on a trajectory  generated 
by the differential equation
(\ref{3.4}), we have
$$
\frac{d}{ds} ~\Re P_j ( {\tilde \chi}) = 
~| P_j^\prime ({\tilde \chi}) | $$
and
hence one of the two conditions in Theorem \ref{T3.1} is satisfied by
the path $\mathcal{C}_j = \left \{ {\tilde \chi} : {\tilde \chi} = \chi_j (t, \chi)
\right \}$, provided it remains within $\mathcal{E}$.
\end{Remark}

\begin{Lemma}
\label{L3.3}
$\Re P_1$ increases monotonically on the boundary of $\mathcal{E}_0$ counterclockwise from
$\infty \e^{i \frac{2}{9} \pi}$ to $\infty \e^{-i \frac{2}{9} \pi}$ with
$$ \frac{d}{ds} ~\Re P_1 ( \chi (s)) ~>~C | \chi (s)|^{5/4} ,~~$$
while 
$\Re P_2$ increases monotonically on the boundary of $\mathcal{E}_0$ clockwise from
$\infty \e^{-i \frac{2}{9} \pi}$ to $\infty \e^{i \frac{2}{9} \pi}$ with
$$ \frac{d}{ds} ~\Re P_2 ( \chi (s)) ~>~C | \chi (s)|^{5/4} $$
$s$ being arc-length on
$\mathcal{E}_0$.
\end{Lemma} 
\begin{proof} Consider the solution to   
(\ref{3.5}), with initial condition on the imaginary $\psi$-axis  slightly  above
 $\psi=i \sqrt{\frac{2}{9}} $. This corresponds to starting at $\chi = \infty
 \e^{i 2 \pi/9} $ with $P_1^\prime (\chi)$ and tracing the Stokes line where
  $\text{{\rm Im} }~P_1=0$ and $\Re P_1$ is increasing. 
  From the
  equation it is clear that $\psi$ remains on the imaginary axis and approaches
  $i \infty$, implying that $\arg \chi = \frac{2}{9} \pi $ 
  is a Stokes line
  where $\Re P_1$ is increasing monotonically all the way to the
  origin in the $\chi$-plane. This also means that locally near $\chi = 0$, $P^\prime \sim 
  \omega_1 \chi^{1/2}$ and
  $P_1 ~\sim ~\frac{2}{3} \omega_1 \chi^{3/2} $, since this is the only root
  of the cubic (\ref{3.3}) which is real on $\chi = r \e^{i 2 \pi/9}$. This corresponds to
  $\psi \sim \omega_1 \chi^{-3/4} $ as $\chi ~\rightarrow ~0$.  Now, 
taking the initial condition slightly above $\psi=i \sqrt{\frac{2}{27}}$, it is clear
  from the differential equation (\ref{3.5}) that $\psi $ 
   remains on the positive imaginary
  $\psi$-axis and approaches $\psi = 
   i \sqrt{\frac{2}{9}} $ from
  below. This corresponds to the fact that $\arg\chi = - \frac{2}{9} \pi$ is a
  Stokes line beyond the turning point $\chi=\chi_s = \left ( \frac{81
      \sqrt{3}}{4\sqrt{2}} \right )^{4/9} \e^{-i 2 \pi/9} $, with $\Re P_1$ increasing monotonically towards $\infty 
\e^{-i \frac{2}{9}
    \pi}$ and for large $r$, $\frac{d}{dr} ~\Re P_1 \ge C r^{5/4} $.  
Now, consider the segment 
$\chi = r ~\e^{-i 2 \pi/9} $, where $0
  ~<~r ~<~ \left ( \frac{81 \sqrt{3}}{4\sqrt{2}} \right )^{4/9} $.  If we
  introduce the transformations
$$ \psi = i\,\, \Psi;~~~~\chi = r ~\e^{-i 2 \pi/9} $$ 
into (\ref{3.2}), then
$$ \Psi^3 - \frac{2}{9} \Psi + q^{-1} = 0 , {\rm where} ~~q = r^{9/4}$$
The roots of the cubic that corresponds to $P_j^\prime$ are:
\begin{equation}
\label{3.6}
\Psi= \Psi_j = -\frac{(2916)^{1/3}}{18 q^{1/3}} J^{1/3} \omega_j^{-1}  
~-~\frac{4 q^{1/3}}{3 (2916)^{1/3} J^{1/3}} \omega_j ~~{\rm where}~~J
= 1 - \sqrt{1 - \frac{96 q^2}{59049}} 
\end{equation}
(the principal branch is used).
The asymptotic behavior of $\Psi_j$ for small $r$ is given by
\begin{equation}
\label{3.6.1}
\Psi_1 \sim \e^{i \pi/3} r^{-3/4}, \Psi_2 \sim -r^{-3/4} 
,\Psi_3 \sim  \e^{-i \pi/3} r^{-3/4}
\end{equation}
From (\ref{3.6}), it follows that on the line $\chi = r \e^{-i 2 \pi/9}$, for
$0 < r <
\left (\frac{81 \sqrt{3}}{4 \sqrt{2}} \right )^{4/9} $,
we have $$\frac{d}{dr} \Re P_1 (r \e^{-i 2\pi/9}) 
= r^{5/4} \Re \Psi_1 >C r^{1/2}> 0$$ 
Thus, for all $r$, we have   
$\frac{d}{dr} {\Re} P_1 \left ( r \e^{-i 2 \pi/9} \right ) > C r^{5/4}$. 
From the reflection-symmetry between 
$P_1$ and $P_2$ on the positive real $\chi$-axis, the
statement for $P_2$ follows. 
\end{proof}
\begin{Lemma}
\label{L3.3}
$\Re P_3$ decreases monotonically on the boundary of 
$\mathcal{E}_0$ counter-clockwise from
$\infty \e^{\pm{i} \frac{2}{9} \pi}$ to 0, and
$$ \frac{d}{ds} \Re P_3 ( \chi (s)) >
\frac{C |\chi(s)|^{1/2}}{|\chi(s)|^{3/2}+1} $$
$s$ being the arc-length towards $\infty$. In this, the positive
real $\chi$-axis is a Stokes line 
with $\Re P_3$ increasing towards $\infty \e^{i 0}$ and
satisfying  the  above monotonicity condition. 
\end{Lemma}
\begin{proof}
Consider (\ref{3.5}) starting with $\psi$ on the positive imaginary
axis,  slightly below $\psi=i \sqrt{\frac{2}{27}}$, corresponding to $\chi =
  \left ( \frac{81\sqrt{3}}{4 \sqrt{2}} \right )^{4/9} \e^{-i 2
    \pi/9}$. The differential equation implies that $\psi$ 
remains on
  the positive imaginary axis 
as it moves towards the origin.  This
corresponds to $\chi = \infty \e^{-i \frac{2}{9} \pi}$, since 
$\psi \sim
  \frac{9}{2} \chi^{-9/4}$ for large $\chi$, where $P_3^\prime \sim
  \frac{9}{2 \chi}$.  Thus, the segment $\chi = r \e^{-i \frac{2}{9}
    \pi}$, $r > \left ( \frac{81\sqrt{3}}{4 \sqrt{2}} \right )^{4/9}$
  is a Stokes line with
$$
\frac{d}{dr} \Re P_3 \left ( r \e^{-i \frac{2}{9} \pi} \right ) > \frac{C}{r} $$
From the symmetry
about the real $\chi$-axis, the same argument can be repeated for $\chi = r
\e^{i \frac{2}{9} \pi}$ for $r > \left ( \frac{81\sqrt{3}}{4 \sqrt{2}} \right
)^{4/9}$ to show that 
this segment is  also part of the Stokes line 
with 
$\Re P_3$ increasing with $r$.

For $r <\left ( \frac{81\sqrt{3}}{4 \sqrt{2}} \right )^{4/9}$, 
an examination of  $\Psi_3$ 
 in (\ref{3.6})
shows that 
$\Re P_3 (r \e^{\pm{i} \frac{2}{9} \pi })$ 
continues to decrease monotonically with decreasing $r$, 
 though these segments are not part of any 
Stokes line. Near the origin, given the asymptotics of $\Psi_3$ 
in (\ref{3.6.1}), 
it follows that $P_3 (\chi) \sim \frac{2}{3} \chi^{3/2} $. Hence a  corresponding
inequality follows,  incorporating this behavior at the origin, while at the same time
satisfying condition for large $\chi$ 
$$ \frac{d}{dr} \Re P_3 \left ( r \e^{\pm{i} \frac{2}{9} \pi} \right ) >
\frac{C r^{1/2}}{r^{3/2}+1} $$
which implies the inequality in the Lemma. 
Furthermore, when $\arg\chi=0$, it is easily seen that
$P_3^\prime$ is real and positive and so $P_3$ increases monotonically to
$\infty$ as we approach $\infty \e^{i 0}$. 
\end{proof}

\begin{Lemma}
\label{L3.4}
  For any $\delta\in (0,\frac{\pi}{63})$ there exists $R_0$ independent of
$\delta$ so that
\begin{enumerate}
\item{}
  $$
  \frac{d}{ds} \Re P_{1} (\chi (s)) \ge C |\chi(s)|^{5/4} $$
  for $C>0$ independent of any parameter, where $s$ is the arc-length
  representation of part of the boundary of $\mathcal{E}$ for which $ |\chi
  (s)| >R_0$;  $s$ is increasing in $r$ when $\chi = r \e^{-i \frac{2}{9} \pi +
    i \delta }$ and decreasing  when $\chi = r \e^{i \frac{2}{9} \pi
    - i\delta }$.
  
\item{} $$
  \frac{d}{ds} \Re P_{2} (\chi (s)) \ge C
  |\chi(s)|^{5/4} $$
  for $C>0$ independent of any parameter, where $s$ is
  the arc-length representation of part of the boundary of $\mathcal{E}$ for
  which $ |\chi (s)| >R_0$;  $s$ is increasing in $r$ when $\chi = r 
\e^{i
    \frac{2}{9} \pi - i \delta }$  and decreasing when $\chi = r
  \e^{-i \frac{2}{9} \pi + i \delta}$.

\item{} $$ \frac{d}{ds} \Re P_{3} (\chi (s)) \ge C |\chi(s)|^{-1} $$ 
for $C>0$ independent of any parameter, where $s$ is the arc-length representation
of part of the boundary of $\mathcal{E}$ for which $ |\chi (s)| >R_0$; $s$
is increasing in $r$ 
when  $\chi =r \exp \left \{ \pm{i} \left [\frac{2}{9} \pi - \delta \right ] \right \}$.
\end{enumerate}
\end{Lemma}

\begin{proof}
  This follows from the asymptotic behavior of $P_1^\prime $,
  $P_2^\prime$ and $P_3^\prime$ for large $\chi$ in (\ref{3.5.1}) after noting
  that
$$ \frac{d}{dr} \Re P_j (r \e^{i \theta} ) = 
\Re \left [ \e^{i \theta} P_j^\prime (r \e^{i \theta} ) \right ]
$$
\end{proof} 

\begin{Lemma} 
\label{L3.5}
For $ 0 < \epsilon_1 \le r \le R_0$. There exists a 
small enough $\delta >0$, independent 
of any parameter, so that 
$$\frac{d}{dr} \Re P_j \left (r \e^{-i \frac{2}{9} \pi + i \delta } 
\right )>C>0
~{\rm for} ~j=1,3
$$ 
while 
$$ -\frac{d}{dr} \Re P_2 \left (r \e^{-i \frac{2}{9} \pi + i \delta } \right )>C>0
$$ 
with $C$ independent of $\delta$. 
Again, 
for $ \epsilon_1 \le r \le R_0$, there is a  $\delta >0$, independent of
any parameter so that  
$$\frac{d}{dr} \Re P_j \left (r \e^{i \frac{2}{9} \pi - i \delta } 
\right )>C>0
~{\rm for}~j=2,3$$ 
while 
$$ -\frac{d}{dr} \Re P_1 \left (r \e^{i \frac{2}{9} \pi - i \delta } 
\right )>C>0$$ 
for some $C$ independent of $\delta$. 
\end{Lemma}

\begin{proof}
  From the lemmas about the behavior of $P_j$ on 
  $\partial \mathcal{E}$, the statements are clearly true for $\delta = 0$. From
  continuity, it follows that the same is true (adjusting $C$) for all
  sufficiently small $\delta $ and hence the lemma follows.
\end{proof}

\begin{Definition}
$$ {\partial \mathcal{E}}_L = {\partial \mathcal{E}}_{L}^+ \cup \mathcal{E}_L^-$$
where 
$$ {\partial \mathcal{E}}_L^+ = \left \{ \chi = \chi_3 + r \e^{i 2\pi/3} ~{\rm for} ~       
 0 \le r \le |\chi_2 - \chi_3| \right \} $$ 
$$ {\partial \mathcal{E}}_L^- = \left \{ \chi = \chi_3 + r \e^{-i 2 \pi/3} ~{\rm for} ~       
 0 \le r \le |\chi_2 - \chi_1 \right \} $$ 
\end{Definition}

\begin{Lemma}
\label{L3.6}
$\Re P_3$ increases in $r$ on  $\partial \mathcal{E}_{L}^+$ and
$\partial \mathcal{E}_L^-$. $\Re P_1 $ decreases in $r$ on $\partial \mathcal{E}_L^+$, but 
increases in $r$ on $\partial \mathcal{E}_L^-$. $\Re P_2$ increases in 
$r$ on $\partial \mathcal{E}_L^+$ and  decreases in $r$ on $\partial
\mathcal{E}_L^-$ and in all cases,  we have  on ${\partial \mathcal{E}}_L$, 
$$ \Big\lvert \frac{d}{dr} \Re P_j (\chi (r) ) \Big\rvert \ge C >0 $$
where $C$ only depends on the choice of $|\chi_j|$. $\Re P_j $ attains a maximum
on ${\partial \mathcal{E}}_L$ at the corresponding $\chi_j$. 
\end{Lemma}
\begin{proof}
We note that since $|\chi_3|$ is small, we have
$$-\frac{d}{dr} \Re P_3 (\chi (r)) =-\Re \left [ P_3^\prime (\chi(r)) \e^{i 2\pi/3} \right ]
\sim |\chi(r)|^{1/2}\sin \left ( \frac{\pi}{6} + \frac{\theta}{2} \right ) >C>0 
$$
where $\arg\chi = \theta \in \left [ - \frac{2 \pi}{9} +\delta,
  \frac{2 \pi}{9} -\delta \right ]$.  By symmetry  we also get for
$\chi $ on ${\partial \mathcal{E}}_L^-$
$$-\frac{d}{dr} \Re P_3 (\chi (r)) =-\Re \left [ P_3^\prime (\chi(r)) 
\e^{-i 2 \pi/3}  \right ]
\sim |\chi(r)|^{1/2}\sin \left ( \frac{\pi}{6} - \frac{\theta}{2} \right ) >C>0 $$
For  $P_1$ 
we find that for $\chi \in {\partial \mathcal{E}}_L^+$,  
$$-\frac{d}{dr} \Re P_1 (\chi (r)) 
\sim |\chi(r)|^{1/2}\cos \left (\frac{\pi}{3} + \frac{\theta}{2} \right ) >C>0 $$
On ${\partial\mathcal{E}}_L^-$, we obtain   
$$\frac{d}{dr} \Re P_1 (\chi (r)) \sim |\chi(r)|^{1/2}\cos \left
(\frac{\theta}{2} \right ) >C>0 $$ Thus, on ${\partial
\mathcal{E}}_L$, $\Re P_1$ increases monotonically from top to bottom
with $\frac{d}{ds} \Re P_1 (\chi (s)) >C>0$. On this boundary $P_2$
increases monotonically from bottom to top by a similar argument. On
the other hand, $P_3$ is maximum at $\chi_3$; it decreases as we move
up or down.
\end{proof}
   
\begin{Lemma} 
\label{L3.7}
On the boundary of $\mathcal{E}$, $\Re P_1 $ increases monotonically with  $s$ as
we traverse the boundary counterclockwise and:
$$
\frac{d}{ds} \Re P_1 (\chi (s)) \ge C |\chi (s)|^{5/4}>0
$$
whereas $ \Re P_2$ increases monotonically with the arclength
$s$ as this boundary is traversed clockwise and
$$ \frac{d}{ds} \Re P_2 (\chi (s)) \ge C |\chi (s)|^{5/4} >0$$
On the other hand at the upper part  of $\partial \mathcal{E}$, i.e.
on  $\partial \left(\mathcal{E} \cap \left \{ 
\chi : \text{{\rm Im} }\chi >0 \right \}\right) $, 
$$ \frac{d}{ds} \Re P_3 (\chi (s)) \ge C |\chi (s)|^{-1} >0$$
where the boundary is traversed counterclockwise. 
For the lower part  of $\partial \mathcal{E}$, i.e. on   $\partial \left(\mathcal{E} \cap \left \{ 
\chi : \text{{\rm Im} }\chi <0 \right \}\right) $ we have
$$ \frac{d}{ds} \Re P_3 (\chi (s)) \ge C |\chi (s)|^{-1} >0$$
where the boundary is now traversed clockwise. 
\end{Lemma}
\begin{proof} The proof follows from   Lemmas \ref{3.4}-\ref{3.6}.
\end{proof}

\noindent{\bf Proof of Theorem \ref{T3.1}.}
Any $\chi \in \partial \mathcal{E}$ can be joined to $\infty \e^{i
  \theta_j}$ along $\partial \mathcal{E}$ so that $\frac{d}{ds} \Re
P_j ({\tilde \chi} (s))$ satisfies the lower bounds given in Lemma \ref{L3.7}.
If $\chi \in \mathcal{E}$, we choose steepest ascent paths for $\Re
P_j$ until (i) it goes to $\infty$, or (ii) it intersects $\partial
\mathcal{E}$, from which point we continue along the ascent paths of
$\partial \mathcal{E}$. The proof is  complete.

\section{Estimates on the solution $G_k$ in the domain $\mathcal{E}_k $}

The main theorem proved in this section is the following.
\begin{Theorem}
\label{T4.1}
For $\eta \in \mathcal{E}_k $ we have 
$$ \| \eta^{3/2} G_k \|_{\infty, \mathcal{E}_k} \le \frac{K}{k} \Big\| 
\frac{R_k}{G_0^3} \Big\|_{\infty, \mathcal{E}_k} $$
$$ \| \eta^{5/2} G_k^\prime \|_{\infty, \mathcal{E}_k} \le K \Big\| 
\frac{R_k}{G_0^3} \Big\|_{\infty, \mathcal{E}_k} $$
$$ \| G_k^{\prime \prime} \|_{\infty, \mathcal{E}_k} \le K \Big\| 
\frac{R_k}{G_0^3} \Big\|_{\infty, \mathcal{E}_k} $$   
where the constant $K$ is independent of $k$ (and therefore of $\beta_k$).
\end{Theorem}

\begin{Remark} The proof comes at the end of \S 5, after a few lemmas. 
  It is convenient to derive an integral equation for $G_k$ and its
  first two derivatives.  We exploit the largeness of $k$ to control
  the solution of the integral equation.  The asymptotic
    behavior of the solution of the homogeneous equation
    $\mathcal{L}_k u = 0$ is different in the regimes $\eta \ll
    k^{4/9} $ and $|\eta| \gtrsim k^{4/9} $.  Consequently,
    different integral equations will be used in $\mathcal{E}_k$ and
    $\mathcal{D}_k $ (analyzed in the next section).
\end{Remark}
In $\mathcal{E}_k$, it is convenient to introduce scaled variables:
\begin{equation}
\label{4.0.0}
\eta = \beta_k^{4/9} \chi ,~G_k (\beta_k^{4/9} \chi ) = z_k (\chi)
\end{equation}
Then, ({\ref{0.2.8}) becomes
\begin{equation}
\label{4.0.1}
\mathcal{\hat L}_k z_k = - \frac{2 \Psi_0}{9 \chi^2} z_k^\prime 
+ \beta_k \left ( \frac{\Psi_0}{\chi^3} + \frac{\Psi_1}{\beta_k
\chi^3} \right ) z_k + \frac{R_k}{G_0^3} \beta_k^{4/3} = \mathcal{R} (\chi), 
\end{equation}
where
\begin{equation}
\label{4.0.2} 
\mathcal{\hat L}_k u :=u^{\prime \prime \prime} + \frac{2}{9} \beta_k^2 \chi^{5/2} u^\prime 
- \beta_k^3 \chi^{3/2} u,  
\end{equation}  
and $\Psi_0$ and  $\Psi_1$ are defined by
\begin{equation}
\label{4.0.2.1}
\frac{1}{G_0^3} - \eta^{3/2} = - \frac{\Psi_0}{\eta^3} ,~
-\frac{3 G_0^2 G_0^{\prime \prime \prime}}{G_0^3} = - \frac{\Psi_1}{\eta^3} 
\end{equation}
From the large $\eta \in \mathcal{E}_k$ behavior of $G_0$ we
  see that $\Psi_0$ and $\Psi_1$ are bounded for large $\beta_k$ as well
  as for large $\chi$.
Let $v$ be the solution for $\chi \in \mathcal{E}$ of
\begin{equation}
\label{4.1} 
\mathcal{\hat L}_k v = \mathcal{R}
\end{equation}
Using rigorous WKB results \cite{Wasow}, it follows that for large $\beta_k$,
there exist three independent solutions of the associated homogeneous
equation, with leading behavior $v_1$, $v_2$, $v_3$ where
\begin{equation}
\label{4.3}
v_j (\chi) = \e^{\beta_k P_j (\chi)  + W_j (\chi)}
\end{equation}
where
$\alpha = P_j^\prime$ are the three roots of the cubic equation
\begin{equation}
\label{4.4}
\alpha^3 + \frac{2}{9} \chi^{5/2} \alpha - \chi^{3/2} = 0
\end{equation}
Note that two roots of (\ref{4.4})  coincide iff $\alpha^2 + \frac{2}{27}
\chi^{5/2} = 0$ i.e. iff
$$\chi =\chi_s = \left ( \frac{81 \sqrt{3}}{4\sqrt{2}} \right )^{4/9}
\e^{\pm{i} 2 \pi/9 } $$
only possible \emph{outside} $\mathcal{E}$. Hence the
$v_i, i=1,2,3$ are independent in $\mathcal{E}$.  The corresponding
$W_j$ are given by
\begin{equation}
\label{4.5}
W_j^\prime = 
-\frac{3 P_j^\prime P_j^{\prime \prime}}{3 P_j^{\prime^2} 
+ \frac{2}{9} \chi^{5/2} } 
\end{equation}
and the $P_j$ are  uniquely determined by the following asymptotic conditions
for 
large $\chi$:  
\begin{equation}
\label{4.6}
P_1 = \frac{4 \sqrt{2}}{27} i \chi^{9/4} - \frac{9}{4} \ln \chi+o(1),
P_2 = -\frac{4 \sqrt{2}}{27} i \chi^{9/4} - \frac{9}{4} \ln \chi+o(1),~
P_3 = \frac{9}{2} \ln \chi +o(1) 
\end{equation}
\begin{equation}
\label{4.6.1}
W_1 = -\frac{15}{8} \ln \chi +o(1)
,~
W_2 = -\frac{15}{8} \ln \chi +o(1)
,~ W_3 =o(1)
\end{equation}
We now use the  $v_i$ to  write an
integral equation  for $v$, equivalent to 
(\ref{4.1}), with appropriate decay conditions at
$\infty$. First, we have 
\begin{equation}
\label{4.7} 
\mathcal{M}:= \left[ \begin {array}{ccc} v_1&v_2&v_3\\
           \beta_k^{-1} v_1^\prime& \beta_k^{-1} v_2^\prime & \beta_k^{-1} v_3^\prime\\
           \beta_k^{-2} v_1^{\prime\prime}&\beta_k^{-2} v_2^{\prime\prime}& 
           \beta_k^{-2} v_3^{\prime\prime}
           \end {array} \right]
= \left [ \begin {array}{ccc} 1 & 1 & 1\\
           m_{21}& m_{22} & m_{23} \\
           m_{31}& m_{32} & m_{33} 
           \end {array} \right ] 
\left[ \begin {array}{ccc} v_1 & 0 & 0\\
           0 & v_2 & 0 \\
           0 & 0 &   v_3 
           \end {array} \right ]
\end{equation}
where for $j=1,2,3$ 
$$ m_{2j} = P_j^\prime + \frac{W_j^\prime}{\beta_k} $$
$$ m_{3j} = \left ( P_j^\prime + \frac{W_j^\prime}{\beta_k} \right )^2 
           + \frac{1}{\beta_k} \left ( P_j^{\prime \prime} + \frac{W_j^{\prime \prime}}{\beta_k} \right ) $$
From the asymptotic properties of $P_j$ and $W_j$,
it follows that for large $\beta_k$ we have
$m_{ij}=O(1)$ for all $i,j$. Furthermore, for large $\chi$, we also have
\begin{equation}
\label{4.7.1}
m_{21} = O(\chi^{5/4} ), m_{22} = O(\chi^{5/4}), m_{23} = O(\chi^{-1})
\end{equation}
\begin{equation}
\label{4.7.2}
m_{31} = O(\chi^{5/2} ), ~m_{32} = O(\chi^{5/2}), ~m_{33} = O(\chi^{-2})
\end{equation}
Let
\begin{equation}
\label{4.8}
Q_1 = ( \mathcal{M}^\prime - Q_2 \mathcal{M}) \mathcal{M}^{-1} 
,~{\rm where} ~
Q_2 =  \beta_k \left[ \begin {array}{ccc} 0 & 1 & 0\\
           0 & 0 & 1 \\
           \chi^{3/2} & -\frac{2}{9} \chi^{5/2} & 0  
           \end {array} \right] 
\end{equation}
Then $\mathcal{M}$ satisfies the differential equation
\begin{equation}
\label{4.9} 
\mathcal{M}^\prime - (Q_2+Q_1) \mathcal{M} = 0 
\end{equation}
Denoting
\begin{equation} 
\label{4.9.0}
\left[ \begin {array}{ccc} n_{11} & n_{12} & n_{13}\\
           n_{21}& n_{22} & n_{23} \\
           n_{31}& n_{32} & n_{33} 
           \end {array} \right ] =  
\left[ \begin {array}{ccc} 1 & 1 & 1\\
           m_{21}& m_{22} & m_{23} \\
           m_{31}& m_{32} & m_{33} 
           \end {array} \right ]^{-1} 
\end{equation}
and 
$$\Delta=m_{22} 
m_{33} - m_{23} m_{32} - m_{21} m_{33} + m_{21} m_{32} 
- m_{31} m_{22} + m_{31} m_{23} $$
we have
\begin{equation}
\label{4.9.0.5}
n_{1,3} = (m_{23}-m_{22})/\Delta;\ \  
n_{2,3} = (m_{21}-m_{23})/\Delta;\ \ n_{3,3} =(m_{21}-m_{22})/\Delta
\end{equation}
The first two rows of $\mathcal{M}^\prime - Q_2 \mathcal{M}$ are zero.
Hence, the 
same is true for the first two rows of $Q_1$. Therefore, 
\begin{equation}
\label{4.9.1}
Q_1 = \beta_k^{-1} \left[ \begin {array}{ccc} 0 & 0 & 0\\
           0 & 0 & 0 \\
           b_{31} & b_{32}  & b_{33}   
           \end {array} \right]
\end{equation}
Defining $r_j$ so that
\begin{equation}
\label{4.10}
\mathcal{\hat L}_k v_j = \beta_k r_j v_j, 
\end{equation}
we have
\begin{equation}
\label{4.10.0}
r_j = 3 P_j^\prime W_j^{\prime^2} + 3 P_j^\prime W_j^{\prime \prime} + 3 W_j^\prime P_j^{\prime \prime}
+ P_j^{\prime \prime \prime} + \beta_k^{-1} \left ( W_j^{\prime^3} + 3 W_j^\prime W_j^{\prime \prime} + 
W_j^{\prime \prime \prime} \right )
\end{equation}
We note that  $r_j = O(1)$ for large $\beta_k$. For large $\chi$ we have
\begin{equation}
\label{4.10.0.1}
r_1 = O(\chi^{-3/4}), ~~r_2 = O(\chi^{-3/4} ) ,~~r_3 = O(\chi^{-5})
\end{equation}
Also, with $\Delta_1={
(m_{21} - m_{22}) (m_{32} - m_{33} ) - (m_{22} - m_{23} ) (m_{31} - m_{32} )}  $we have
\begin{equation} 
\label{4.10.1}
b_{3,2} = [(r_1 - r_2) (m_{32} - m_{33}) - (r_2 - r_3) (m_{31} - m_{32})]/\Delta_1
\end{equation}
\begin{equation}
\label{4.10.2}
b_{3,3} = -[(r_1 - r_2) (m_{22} - m_{23}) - (r_2 - r_3) (m_{21} - m_{22})]/\Delta_1
\end{equation}
\begin{equation}
\label{4.10.3}
b_{3,1} = r_3 - b_{3,2} m_{2,3} - b_{3,3} m_{3,3}
\end{equation}
From the asymptotics of $r_j$ and $m_{i,j}$ for large $\beta_k$ we get 
$b_{3,j} = O(1)$. For large 
$\chi \in \mathcal{E}$ we have
\begin{equation}
\label{4.10.4}
b_{3,2} = O(\chi^{-2}), ~b_{3,3} =  O(\chi^{-2}) ,~~{\rm and} ~~b_{3,1}= O(\chi^{-3})
\end{equation}
Then,  for large $\chi \in \mathcal{E}$, it follows that
\begin{equation}
\label{4.10.5}
n_{1,3} = O(\chi^{-5/2}), ~~n_{2,3} = O(\chi^{-5/2}) ~~{\rm and} ~~n_{3,3} = O(\chi^{-5/2})
\end{equation}
In the domain $\mathcal{E}$ the $b_{3,j}$ are
analytic in $\chi$, bounded for large $\beta_k$ 
and decay for large $\chi$, (see (\ref{4.10.4})).
Furthermore, (\ref{4.9}) implies 
\begin{equation}
\label{4.10.6}
\mathcal{L}_{_{WKB}} v_j := v_j^{\prime \prime\prime} - \frac{b_{33}}{\beta_k} v_j^{\prime \prime} 
+ \left [ \frac{2}{9} \beta_k^2 \chi^{5/2} - b_{3,2} \right ] v_j^\prime  
- \left [ \beta_k^3 \chi^{3/2} + \beta_k b_{3,1} \right ] v_j = 0 
\end{equation}
Using variation of parameters, we see 
that  one solution of $\mathcal{\tilde
  L}_k v = \mathcal{R}$  satisfies:
\begin{equation}    
\label{4.11}
v (\chi) = \mathcal{\tilde V} \left [ {\hat R} \right ] (\chi)  ~~;~~
v^\prime (\chi) = \mathcal{\tilde V}^\prime \left [ {\hat R} \right ] (\chi) ~~;~~ 
v^{\prime\prime} (\chi) = \mathcal{\tilde V}^{\prime\prime} \left [  {\hat R} \right ] (\chi) ~~;~~ 
\end{equation}
where
\begin{equation}
\label{4.12}
{\hat R} (\chi) = \mathcal{R} (\chi) -\beta_k^{-1} b_{3,3} v^{\prime \prime} 
- b_{3,2} v^\prime (\chi) -\beta_k b_{3,1} v (\chi) ,~~ 
\end{equation}
and the operators $\mathcal{\tilde V}$, $\mathcal{\tilde V}^\prime$ and
$\mathcal{\tilde V}^{\prime \prime}$ are defined by:  

\begin{equation}
\label{4.13}
\mathcal{\tilde V} [\hat R] (\chi) = 
\sum_{j=1}^3 \frac{1}{\beta_k^2} \int_{\infty \e^{i \theta_j}}^\chi  
\e^{\beta_k ~\left [P_j (\chi) - P_j ({\tilde \chi}) \right ] + W_j (\chi)
~-~W_j ({\tilde \chi}) } n_{j,3} ({\tilde \chi}) {\hat R} ({\tilde \chi}) d{\tilde \chi}
\end{equation}
\begin{equation}
\label{4.13.1}
\mathcal{\tilde V}^{\prime} [\hat R] (\chi) =
\sum_{j=1}^3 \frac{m_{2,j} (\chi)}{\beta_k} \int_{\infty \e^{i \theta_j}}^\chi  
\e^{\beta_k ~\left [P_j (\chi) - P_j ({\tilde \chi}) \right ] + W_j (\chi)
~-~W_j ({\tilde \chi}) } n_{j,3} ({\tilde \chi}) {\hat R} ({\tilde \chi}) d{\tilde \chi}
\end{equation}
\begin{equation}
\label{4.13.2}
\mathcal{\tilde V}^{\prime \prime} [\hat R] (\chi) =
\sum_{j=1}^3 m_{3,j} (\chi)  
\int_{\infty \e^{i \theta_j}}^\chi  
\e^{\beta_k ~\left [P_j (\chi) - P_j ({\tilde \chi}) \right ] + W_j (\chi)
~-~W_j ({\tilde \chi}) } n_{j,3} ({\tilde \chi}) {\hat R} ({\tilde \chi}) d{\tilde \chi}
\end{equation}
where $\theta_1 = -\frac{2}{9} \pi + \delta $, $\theta_2 = \frac{2}{9} \pi -
\delta $ and $\theta_3 = 0 $, and the paths of integration $\mathcal{C}_j$ are
chosen to be the ascent paths for $\Re P_j$  of   Theorem \ref{T3.1}.
Also, note that for large $\chi$, $W_j (\chi)$ grows at most logarithmically
with $\chi$ implying that, uniformly in $\mathcal{E}$, we have $ W_j
=o(\beta_k P_j) $. As we shall see, there is a unique solution $v$ of
  (\ref{4.11}) that decays as $\chi \rightarrow \infty$ in $\mathcal{E}$, with
  $\hat{R}$ 
  having similar decay properties. The way we show
  this is by proving contractivity of the integral system in a suitable space
  of decaying functions.  (In fact, there can be no other decaying solutions,
  since the associated homogeneous equation does not have nonzero decaying
  solutions in $ \mathcal{E}$.)
\begin{Lemma}
\label{L4.1}
If the $P_j $ satisfy Property 1 in $\mathcal{E}$, then for sufficiently large
$k$ (or, which amounts to the same, large $\beta_k$)  we have 
\begin{equation}
\label{4.16}
\| \chi^{3/2} \mathcal{\hat V} [ \hat R ] \|_{\infty, \mathcal{E}} \le ~
\frac{C}{\beta_k^3} ~\| \hat R \|_\infty 
\end{equation}
\begin{equation}
\label{4.17}
\| \chi^{5/2} \mathcal{\hat V}^\prime [ \hat R ] \|_{\infty, \mathcal{E}}
\le ~\frac{C}{\beta_k^2} ~\| \hat R \|_\infty 
\end{equation}
\begin{equation}
\label{4.17.1}
\| \mathcal{\hat V}^{\prime \prime} [ \hat R ] \|_{\infty, \mathcal{E}}
\le ~\frac{C}{\beta_k} ~\| \hat R \|_\infty 
\end{equation}
where the constant $C$ is independent of $\hat R$ and $\beta_k$.
\end{Lemma}

\begin{proof} Theorem \ref{T3.1} shows that on  $\mathcal{C}_j$ 
(defined before Lemma~\ref{L4.1}) we have $|{\tilde \chi}|> C |\chi|$ and
$$
\frac{d}{ds} \Re P_{1,2} ({\tilde \chi} (s)) >C |{\tilde \chi} |^{5/4},~
\frac{d}{ds} \Re P_3 ({\tilde \chi} (s)) >C |{\tilde \chi} (s) |^{-1} $$
 Since $W_j^\prime /P_j^\prime$ is bounded, 
this implies that for sufficiently large $\beta_k$ we have
$$  
\frac{d}{ds} \Re \left [ P_{1,2} + \frac{W_{1,2}}{\beta_k} \right ]({\tilde \chi} (s))  >
\frac{C}{2} |{\tilde \chi} |^{5/4} $$
$$  
\frac{d}{ds} \Re \left [ P_{3} + \frac{W_{3}}{\beta_k} \right ] ({\tilde \chi} (s))  >
\frac{C}{2} |{\tilde \chi (s)} |^{-1} $$
Also,
from (\ref{4.7.1}) and (\ref{4.7.2}),
$$| m_{2,1} | ~<~ C ~|\chi|^{5/4},  |m_{2,2} | ~<~ C ~|\chi|^{5/4},~ |m_{2,3} | ~<~ C ~|\chi|^{-1} ,~~$$
$$| m_{3,1} | ~<~ C ~|\chi|^{5/2} ,~| m_{3,2} | ~<~ C ~|\chi|^{5/2},~ |m_{3,3} | ~<~ C ~|\chi|^{-2},$$
while from (\ref{4.10.5}), $ | n_{3,j}| < C |\chi|^{-5/2} $.  Then, 
\begin{multline}
  \frac{1}{\beta_k^2} \left \lvert \int_{\infty \e^{i \theta_j}}^\chi \exp
    \left [ \beta_k ( P_j (\chi) - P_j ({\tilde \chi} ) \right ] n_{3,j}
    R ({\tilde \chi}) \e^{W_j (\chi)-W_j ({\tilde \chi})} d{\tilde \chi} \right\rvert \\
  \le \frac{C \| \hat R \|_{\infty} |\chi|^{-3/2}}{\beta_k^3} \int_0^1 d \left
    [ \exp \left ( \beta_k [ \Re  P_j(\chi)-\Re P_j({\tilde \chi}) ] \right ) \right ]
\end{multline} 
The bounds for $\mathcal{\hat V}$ follow;  $\mathcal{\hat V}^\prime$ and
$\mathcal{V}^{\prime \prime} $ are bounded similarly.
\end{proof}

\begin{Corollary}
\label{cor31}
Define the operator $\mathcal{T}_k $ acting on triples
$\left (z_k, z_k^\prime, z_k^{\prime \prime} \right )$ as follows:   
\begin{equation}
\label{4.18}
\mathcal{T}_k \left ( z_k, z_k^\prime, z_k^{\prime \prime} \right ) (\chi) =
- \beta_k^{-1} b_{3,3} z_k^{\prime \prime} - 
\left ( \frac{2 \Psi_0}{9 \chi^2} + b_{3,2} \right ) z_k^\prime + \beta_k 
\left ( \frac{\Psi_0}{\chi^3} + \frac{\Psi_1}{\beta_k \chi^3} - b_{3,1} \right ) z_k 
\end{equation}
Then, it follows 
$$ 
\| \chi^{3/2} \mathcal{\hat V} \left [ 
\mathcal{T}_k \left ( z_k, z_k^\prime, z_k^{\prime \prime} \right ) \right ] \|_{\infty, \mathcal{E}} 
\le C \left [ \beta_k^{-4} \| z_k^{\prime \prime} \|_{\infty, \mathcal{E}}  
+ \beta_k^{-3} \| \chi^{5/2} z_k^{\prime} \|_{\infty, \mathcal{E}} 
+ \beta_k^{-2} \| \chi^{3/2} z_k \|_{\infty, \mathcal{E}} \right ]
$$  
$$ 
\| \chi^{5/2} \mathcal{\hat V}^\prime \left [ \mathcal{T}_k 
\left ( z_k, z_k^\prime, z_k^{\prime \prime}
\right ) \right ] \|_{\infty, \mathcal{E}} 
\le C \left [ \beta_k^{-3} \| z_k^{\prime \prime} \|_{\infty, \mathcal{E}}  
+ \beta_k^{-2} \| \chi^{5/2} z_k^{\prime} \| + \beta_k^{-1} \| \chi^{3/2} z_k \|_{\infty, \mathcal{E}} \right ]
$$  
$$ 
\| \mathcal{\hat V}^{\prime\prime} \left [ T_k \left ( z_k, z_k^\prime, z_k^{\prime \prime} 
\right ) \right ] \|_{\infty, \mathcal{E}} 
\le C \left [ \beta_k^{-2} \| z_k^{\prime \prime} \|_{\infty, \mathcal{E}}  
+ \beta_k^{-1} \| \chi^{5/2} z_k^{\prime} \| + \| \chi^{3/2} z_k \|_{\infty, \mathcal{E}} \right ]
$$  
\end{Corollary}
\begin{proof} This follows from  Lemma \ref{L4.1} and bounds on
$b_{3,j}$ in (\ref{4.10.4}) and those on $\Psi_0$, $\Psi_1$ that follow from (\ref{4.0.2.1}). 
\end{proof}
\begin{Lemma}
\label{L4.2}
$$\left\| \chi^{3/2} \mathcal{\hat V} \left [ \beta_k^{4/3} \frac{R_k}{G_0^3} (\beta_k^{4/9} \chi) \right ] \right\|_{\infty, \mathcal{E}} 
\le \frac{C\beta_k^{4/3}}{\beta_k^3} \left\| \frac{R_k}{G_0^3} \right\|_{\infty, \mathcal{E}} $$
$$\left\| \chi^{5/2} \mathcal{\hat V}^\prime  
\left [ \beta_k^{4/3} \frac{R_k}{G_0^3} (\beta_k^{4/9} \chi) \right ] \right\|_{\infty, \mathcal{E}} 
\le \frac{C\beta_k^{4/3}}{\beta_k^2} \left\| \frac{R_k}{G_0^3} \right\|_{\infty, \mathcal{E}} $$
$$\left\| \mathcal{\hat V}^{\prime \prime} 
\left [ \beta_k^{4/3} \frac{R_k}{G_0^3} (\beta_k^{4/9} \chi) \right ] \right\|_{\infty, \mathcal{E}} 
\le \frac{C\beta_k^{4/3}}{\beta_k} \left\| \frac{R_k}{G_0^3} \right\|_{\infty, \mathcal{E}} $$
\end{Lemma}

\begin{proof} This is a consequence of Lemma \ref{L4.1}, with ${\hat R}$ replaced by
$\beta_k^{4/3} {R_k}/{G_0^3} $.
\end{proof}

\begin{Lemma}
\label{L4.3}
For $\left\| {R}_k G_0^{-3} \right\|_{\infty, \mathcal{E}} <\infty$, and for 
$\beta_k$ sufficiently large,  
the  system (\ref{4.11}) has a unique solution 
$\left ( z_k (\chi), z^{\prime} (\chi), z^{\prime \prime} (\chi) \right )$ 
in $\mathcal{E}$,
which satisfies the bounds 
$$\| \chi^{3/2} z_k \|_{\infty, \mathcal{E}} 
\le \frac{C\beta_k^{4/3}}{\beta_k^3} \left\| \frac{R_k}{G_0^3} \right\|_{\infty, \mathcal{E}} $$
$$\| \chi^{5/2} z_k^\prime \|_{\infty, \mathcal{E}} 
\le \frac{C\beta_k^{4/3}}{\beta_k^2} \left\| \frac{R_k}{G_0^3} \right\|_{\infty, \mathcal{E}} $$
$$\| z_k^{\prime \prime} \|_{\infty, \mathcal{E}} 
\le \frac{C\beta_k^{4/3}}{\beta_k} \left\| \frac{R_k}{G_0^3} \right\|_{\infty, \mathcal{E}} $$
\end{Lemma}

\begin{proof}
  Define the Banach space $\mathcal{F}$ of triples of functions
    $(u,v,w)$ analytic in the interior of $\mathcal{E}$ and
    continuous in its closure in the norm
  $$
  \| (u, v, w) \|_{\mathcal{F}} = \beta_k^{5/3} \| \chi^{3/2} u
  \|_{\infty, \mathcal{E}} + \beta_k^{2/3} \| \chi^{5/2} v \|_{\infty,
    \mathcal{E}} + \beta_k^{-1/3} \| w \|_{\infty, \mathcal{E}}
 $$
  We associate $z_k$, $z_k^\prime$ and $z_k^{\prime\prime}$
  with $u$, $v$ and $w$ respectively, and consider ${\hat R}$ as
  depending on $u$, $v$ and $w$ for fixed $R_k/G_0^3$.  We define the
  linear operator $\mathbf{L} : \mathcal{F} \rightarrow \mathcal{F}$
  by
$$ \mathbf{L} \left [ (u,v,w) \right ] := \left ( \mathcal{\hat V} \left
[{\hat R} (u,v,w) \right ], \mathcal{\hat V}^\prime \left [{\hat R} (u,v,w)
\right ], \mathcal{\hat V}^{\prime\prime} \left [{\hat R} (u,v,w) \right ]
\right ) $$ where $\mathcal{\hat V}, \mathcal{\hat V}^\prime, \mathcal{\hat
V}^{\prime \prime}$ are now thought of as acting on $u=z_k$, $v=z_k^\prime$
$w=z_k^{\prime \prime}$ for fixed $R_kG_0^{-3}$.  From (\ref{4.0.1}),
(\ref{4.11}) and (\ref{4.12}), the definition of $\mathcal{T}_k$ in
(\ref{4.18}), and the estimates in Corollary \ref{cor31}, it is easily seen
that
\begin{multline*}
\| \mathbf{L} \left [ (u,v,w) \right ]  
 - \mathbf{L} \left [ ({\tilde u}, {\tilde v}, {\tilde w}) \right ] \|_{\mathcal{F}} \\
= \| \mathcal{V} \left [ \mathcal{T}_k (u-{\tilde u}, v-{\tilde v}, w-{\tilde w} ) \right ], 
\mathcal{\hat V}^\prime \left [\mathcal{T}_k (u-{\tilde u},v-{\tilde v},w-{\tilde w}) \right ], 
\mathcal{\hat V}^{\prime\prime} \left [\mathcal{T}_k (u-{\tilde u},v-{\tilde v},w-{\tilde w}) \right ] 
\|  \\
\le \frac{C}{\beta_k^2}   
\| \left ( u-{\tilde u}, v-{\tilde v}, w-{\tilde w} \right ) \|_{\mathcal{F}} 
\end{multline*}
Hence $\mathbf{L}$ is contractive and the system (\ref{4.11}) has a unique solution
$(z_k, z_k^\prime, z_k^{\prime \prime})$. The estimates on $z_k$, $z_k^\prime$, $z_k^{\prime \prime}$ 
follow easily from Lemma \ref{L4.2}.
\end{proof}

\noindent{\bf Proof of Theorem \ref{T4.1}}. 
This is a consequence of Lemma \ref{L4.3}, noting that
$$\eta^{3/2} G_k (\eta) = \beta_k^{2/3} \chi^{3/2} z_k (\chi) , ~\eta^{5/2}
G_k^\prime (\eta) = \beta_k^{2/3} \chi^{5/2} z_k^\prime (\chi), ~ G_k^{\prime
  \prime} (\eta) = \beta_k^{-8/9} z_k^{\prime \prime} (\chi)$$

\section{Estimate of $G_k$ for large $k$ in $\mathcal{D}_k$}

In this section we prove the following.
                    
\begin{Theorem}
\label{T5.1}
In  $\mathcal{D}_k$ (see Definition~\ref{defDk}) we have
\begin{equation}
\label{5.0.3} 
\| \eta^{3/2} G_k (\eta) \|_{\infty, \mathcal{D}_k } \le \frac{K_{10}}{k} \left\|\frac{R_k}{G_0^3} \right\|_{\infty, \mathcal{D}_k \cup \mathcal{E}_k }
\end{equation}
\begin{equation}
\label{5.0.4} 
\| \eta^{5/2} G_k^\prime (\eta) \|_{\infty, \mathcal{D}_k} \le K_{11}
\epsilon^{3/2} 
\left\|\frac{R_k}{G_0^3} \right\|_{\infty, \mathcal{D}_k \cup \mathcal{E}_k }
\end{equation}
\end{Theorem}

\begin{Remark} 
  The proof of theorem (\ref{T5.1}) is completed at the end of \S6, after a few
  lemmas establishing the properties of
  $\mathcal{L}_k^{-1} $.
\end{Remark}
We first find  a representation of the solution to 
\begin{equation}
\label{5.1} 
\mathcal{L}_k G_k = G_k^{\prime \prime \prime} 
+ \frac{2}{9 G_0^3} \eta G_k^\prime - \frac{7 k -1}{9 G_0^3} G_k
+ \frac{3 G_0^{\prime \prime \prime}}{G_0} G_k  = \frac{R_k}{G_0^3}
\end{equation}
for large $k$
for $\eta \in \mathcal{D}_k$ 
where 
$|\eta| $ is  small  compared to $\beta_k^{4/9}$.
Again following  \cite{Wasow},  there exist
three independent solutions $u_1$, $u_2$, $u_3$ 
to the homogeneous equation
$\mathcal{L}_k u = 0$ such that, for
large $\beta_k$ we have
\begin{multline}
\label{5.3}
u_j (\eta) \sim g_j (\eta) = G_0 (\eta) \e^{\omega_j \beta_k^{1/3} P (\eta) },~{\rm where}~ \omega_1 = \e^{i 2 \pi/3},~\omega_2 = \e^{-i 2 \pi/3},~
\omega_3 = 1 ~{\rm and} \\ 
P(\eta) = \int_{\eta_0}^{\eta} \frac{1}{G_0(\eta')} d\eta' ~\text{for fixed }
\eta_0 \in\mathcal{D}_k  
\end{multline}
We use $g_1$, $g_2$, $g_3$ to find a suitable
 integral equation  for the solution $u$ to 
(\ref{5.1}). As in  \S5, it is convenient  to define 
\begin{equation}
\label{5.4} 
\mathcal{M}: = \left[ \begin {array}{ccc} g_1&g_2&g_3\\
           \beta_k^{-1/3} g_1^\prime& \beta_k^{-1/3} g_2^\prime & \beta_k^{-1/3} g_3^\prime\\
           \beta_k^{-2/3} g_1^{\prime\prime}&\beta_k^{-2/3} g_2^{\prime\prime}& \beta_k^{-2/3} g_3^{\prime\prime}
           \end {array} \right]
\end{equation}
and 
\begin{equation}
\label{5.5}
Q_1: = ( \mathcal{M}^\prime - Q_2 \mathcal{M}) \mathcal{M}^{-1} 
,~{\rm where} ~
Q_2: =  \beta_k^{1/3} \left[ \begin {array}{ccc} 0 & 1 & 0\\
           0 & 0 & 1 \\
           \frac{1}{9 G_0^3} - \frac{3 G_0^{\prime\prime\prime}}{\beta_k G_0} & -\frac{2 \eta}{9 \beta_k^{2/3} G_0^3} & 0  
           \end {array} \right] 
\end{equation}
We get 
\begin{equation}
\label{5.6} 
\mathcal{M}^\prime - (Q_2+Q_1) \mathcal{M} = 0 
\end{equation}
Using (\ref{5.3}) we see that 
\begin{equation}  
\label{5.7}
Q_1 = \beta_k^{-2/3} \left[ \begin {array}{ccc} 0 & 0 & 0\\
           0 & 0 & 0 \\
           \frac{G_0^{\prime^3}}{G_0^3} - 2 \frac{G_0^\prime G_0^{\prime\prime}}{G_0^2} + 4 \frac{G_0^{\prime \prime \prime}}{G_0}
& 
           \left [ -\frac{G_0^{\prime^2}}{G_0^2} + \frac{2 G_0^{\prime\prime}}{G_0} + 
\frac{2 \eta}{9 G_0^3} \right ] \beta_k^{1/3} & 0 
           \end {array} \right]
\end{equation}
The columns of $\mathcal{M}$ also satisfy (\ref{5.6}); it
  follows that for j=1,2,3 we have 
\begin{equation}
\label{5.8}
\mathcal{\hat L}_k g_j := g_j^{\prime \prime \prime} + \left ( \frac{G_0^{\prime^2}}{G_0^2} - \frac{2 G_0^{\prime\prime}}{G_0}
\right ) g_j^\prime + \left ( - \frac{\beta_k}{G_0^3} - \frac{G_0^{\prime^3}}{G_0^3}-  \frac{G_0^{\prime\prime\prime}}{G_0}
+ \frac{2 G_0^{\prime} G_0^{\prime \prime}}{G_0^2} \right ) g_j = 0 
\end{equation}
We note that
\begin{equation}
\label{5.9}
\mathcal{\hat L}_k u = \mathcal{L}_k u -b_{3,2} u^\prime - b_{3,1} u
\end{equation}
where
$$
b_{3,2} = \left ( \frac{2 G_0^{\prime\prime}}{G_0} + \frac{2 \eta}{9 G_0^3} 
-  \frac{G_0^{\prime^2}}{G_0^2} \right ) 
$$ 
$$
b_{3,1} = \left (\frac{G_0^{\prime^3}}{G_0^3} +  4 \frac{G_0^{\prime\prime\prime}}{G_0}
- \frac{2 G_0^{\prime} G_0^{\prime \prime}}{G_0^2} \right ) 
$$ 
For large $|\eta|$ in $\mathcal{D}_k$ we find
\begin{equation}
\label{5.9.1}
b_{3,2} = O(\eta^{5/2}) ,~b_{3,1} = O(\eta^{-3})
\end{equation}
Also, $b_{3,1}$ and $b_{3,2}$ are analytic in $\mathcal{D}_k$.  
It follows that the  $G_k $ in (\ref{5.1}) also satisfy the
integral equation
\begin{equation}
\label{5.10}
G_k (\eta) = \mathcal{V} \left [ {\hat R}_k \right ] (\eta) 
+ \sum_{j=1}^3 a_j g_j (\eta)
\end{equation}
where 
\begin{equation}
\label{5.11}
{\hat R}_k (\eta) = \frac{R_k}{G_0^3} - b_{3,2} G_k^\prime - b_{3,1} G_k 
\end{equation}
The constants $a_j$ are defined in (\ref{5.13}) in terms of
$G_k (\eta_{1,k})$, $G_k(\eta_{2,k})$ and $G_k (\eta_{3,k})$
and the operator $\mathcal{V}$ is defined by    
\begin{equation}
\label{5.12}
\mathcal{V} [{\hat R}_k] (\eta) = 
\sum_{j=1}^3 \frac{\beta_k^{-2/3}}{3} \omega_j G_0(\eta) 
\int_{\eta_j}^\eta G_0(\eta') {\hat R}_k(\eta') 
\e^{\omega_j \beta_k^{1/3} [P(\eta) - P(\eta')]} d\eta'    
\end{equation}
The contours of integration chosen in (\ref{5.12}) are ascent
  paths of $\Re \left [ \omega_j P \right ]$, see Corollary
  \ref{Cor:Property1}).
Given  $G_k(\eta_{1,k})$, $G_k(\eta_{2,k})$ and
$G_k(\eta_{3,k})$ we define $a_1$, $a_2$, $a_3$ by
\begin{equation}
\label{5.13}
\left[ \begin {array}{ccc} 
           g_1 (\eta_{1,k}) & g_2 (\eta_{1,k}) & g_3 (\eta_{1,k}) \\
           g_1 (\eta_{2,k}) & g_2 (\eta_{2,k}) & g_3 (\eta_{2,k}) \\
           g_1 (\eta_{3,k}) & g_2 (\eta_{3,k}) & g_3 (\eta_{3,k}) 
           \end {array} \right] 
\left[ \begin {array}{c} 
           a_1 \\
           a_2 \\
           a_3  
           \end {array} \right] 
=\left[ \begin {array}{c} 
           G_k(\eta_{1,k}) - \mathcal{V} [{\hat R}] [\eta_{1,k}] \\
           G_k (\eta_{2,k}) - \mathcal{V} [{\hat R}] [\eta_{2,k}] \\
           G_k (\eta_{3,k}) - \mathcal{V} [{\hat R}] [\eta_{3,k}] 
           \end {array} \right] 
\end{equation}
Again, 
\begin{equation}
\label{5.14}
G_k^\prime (\eta) = \mathcal{V}^\prime [{\hat R}_k ] (\eta)
                   + \beta_k^{1/3} \sum_{j=1}^3 a_j h_j g_j (\eta) 
\end{equation}
where
\begin{equation}
\label{5.15}
\mathcal{V}^\prime \left [ R \right ] (\eta)
= \sum_{j=1}^3 \frac{\beta_k^{-1/3}}{3} \omega_j h_j (\eta) 
G_0(\eta) \int_{\eta_{j,k}}^\eta G_0(\eta') R (\eta') 
\e^{\omega_j \beta_k^{1/3} [P(\eta) - P(\eta')]} d\eta' ,  
\end{equation}  
and
\begin{equation}
\label{5.16}
h_j (\eta) = \frac{\omega_j}{G_0} + \frac{G_0^\prime}{\beta_k^{1/3} G_0} 
\end{equation}
It is to be noted that 
$$ |\eta^{-1/2} h_j (\eta) | <C $$
for some constant $C$ independent of $\beta_k$.

A few properties of $\mathcal{V}$ and $\mathcal{V}^\prime $ follow from
Property 1 of $P(\eta)$ (established in \S3}).

\begin{Lemma}
\label{L5.2}
Assume $ \| R \|_{\infty, \mathcal{D}_k} < \infty$. Then, 
$$
\| \eta^{3/2} \mathcal{V} \left [ R \right ] (\eta) \|_\infty \le
\frac{K_1}{\beta_k} \| R \|_\infty $$
for a constant $K_1$ independent of
$\beta_k$.
\end{Lemma}
\begin{proof}
 Note that on any of the contours $\mathcal{C}_j$,  from Property 1,
 there exists a constant $C >0$ so that  $\tilde{\eta}>C|\eta|$ for
  $\tilde{\eta}\in\mathcal{C}_j$ and 
$$ \frac{d}{ds} \Re \left \{ \omega_j P ({\tilde \eta} (s)) \right \} 
>C_1 |{\eta} (s) |^{1/2} >0 $$
where $s$ is the arc length.
Therefore the proof follows from the estimate
\begin{multline*}
\left\lvert \beta_k^{-2/3} G_0 (\eta) 
\int_{\eta_{j,k}}^\eta G_0 (\eta') R(\eta') \e^{\beta_k^{1/3} (\omega_j [P(\eta)
-P(\eta')])} 
d \eta' \right \rvert \\
\le \int_0^1 d \left [ 
\exp [\beta_k^{1/3} (\Re (\omega_j [P(\eta) - P(\eta')]) ] \right ] 
\frac{C}{\beta_k |\eta|^{3/2}}
\| R \|_{\infty,\mathcal{D}_k} 
\end{multline*}
\end{proof}

\begin{Lemma} 
\label{L5.3} 
Assume $ \| R \|_{\infty, \mathcal{D}_k} <\infty$. Then
$$ \| \eta^{5/2} \mathcal{V}^\prime \left [ R \right ] (\eta) \|_\infty
\le \frac{K_2 |\eta_{3,k}|^{3/2}}{\beta_k^{2/3}} \| R \|_{\infty, \mathcal{D}_k}  $$
where $K_2$ is a constant independent of $\beta_k$.
\end{Lemma}

\begin{proof}
As before, there
exists a constant $C >0$ so that on the contour $\mathcal{C}_j$ we have 
$ C|\eta| <\eta'$ and
$$ \frac{d}{ds} \Re \omega_j P (\eta' (s)) >C |\eta' (s) |^{1/2} >0 $$
where $s$ is the arc length. 
Thus
\begin{multline*}
\left\lvert \beta_k^{-1/3} h_j (\eta) G_0 (\eta) \int_{\eta_{j,k}}^\eta G_0
(\eta') R(\eta') \e^{\omega_j \beta_k^{1/3} P(\eta)-P(\eta')} d \eta'
\right\rvert \\ \le \int_0^1 d \left \{ \exp [\beta_k \Re (P(\eta) - P(\eta'))
] \right \} \frac{C |\eta|^{3/2} }{\beta_k^{2/3} |\eta|^{5/2}} \| R
\|_{\infty,\mathcal{D}_k}
\end{multline*}
The Lemma follows by noting that in $\mathcal{D}_k$ we have 
$|\eta| \le |\eta_{3,k}|$. \end{proof}

\begin{Corollary}  We have 
\label{L5.4} 
$$ \| \eta^{3/2} \mathcal{V} \left [ b_{3,2} G_k^\prime + b_{3,1} G_k \right ] (\eta)
\|_{\infty,\mathcal{D}_k} \le \frac{K_3}{\beta_k} \left [ \| \eta^{5/2} G_k^\prime \|_{\infty, \mathcal{D}_k}
+ \| \eta^{3/2} G_k \|_{\infty,\mathcal{D}_k} \right ]
$$
\end{Corollary}

\begin{proof} This follows from Lemma \ref{L5.2}, and the bounds on $b_{3,2}$ and $b_{3,1}$ in 
(\ref{5.9.1}). \end{proof}

\begin{Corollary} We have 
\label{L5.5} 
$$
\| \eta^{5/2} \mathcal{V}^\prime \left [ b_{3,2} G_k^\prime + b_{3,1}
  G_k \right ] (\eta) \|_{\infty,\mathcal{D}_k} \le K_ 4
\frac{|\eta_{3,k}|^{3/2}}{\beta_k^{2/3}} \left [ \| \eta^{5/2} G_k^\prime
  \|_{\infty,\mathcal{D}_k} + \| \eta^{3/2} G_k \|_\infty \right ]
$$
\end{Corollary}

\begin{proof} This follows from Lemma \ref{L5.3}, and the bounds on $b_{3,2}$ and $b_{3,1}$ in (\ref{5.9.1}).
\end{proof}

\begin{Corollary} The following inequality holds
\label{L5.6} 
$$\left \| \eta^{3/2} \mathcal{V} \left [ \frac{R_k}{G_0^3} \right ] (\eta) \right \|_{\infty, \mathcal{D}_k}
\le \frac{K_5}{k} \left\| \frac{R_k}{G_0^3} \right \|_{\infty, \mathcal{D}_k} $$
for a constant $K_5$ independent of $k$. 
\end{Corollary}
 
\begin{proof} This follows from Lemma \ref{L5.2}.
\end{proof}

\begin{Corollary} We have 
\label{L5.7} 
$$ \left \| \eta^{5/2} \mathcal{V}^\prime \left [ \frac{R_k}{G_0^3} \right ] (\eta) \right \|_{\infty, \mathcal{D}_k}
\le \frac{K_5 \eta_{3,k}^{3/2}}{\beta_k^{2/3}} \left\| \frac{R_k}{G_0^3} \right \|_{\infty, \mathcal{D}_k}  $$ 
\end{Corollary}
 
\begin{proof} This follows from Lemma \ref{L5.3},  after noting that for $\eta \in \mathcal{D}_k$,
$|\eta| \le \eta_{3,k}$. 
\end{proof}

\begin{Definition}
\label{D5.7.0}
Define the linear operators $\mathcal{T}_1$ and $\mathcal{T}_2$ by
$$ \mathcal{T}_1 \left [ G_k, G_k^\prime \right ] (\eta) =
\sum_{j=1}^3 a_j g_j (\eta) $$
$$ \mathcal{T}_2 \left [ G_k, G_k^\prime \right ] (\eta) = \beta_k^{1/3}
\sum_{j=1}^3 a_j h_j (\eta) g_j (\eta) $$ (see (\ref{5.13}))  since  
$\eta_{j,k} \in \partial \mathcal{E}_k $, $G_k (\eta_{j,k})$ are  known
from the previous section.

\end{Definition}

\begin{Lemma}  We have 
\label{L5.8}
$$ \| \eta^{3/2} \mathcal{T}_1 \left [ G_k, G_k^\prime \right ] \|_{\infty, \mathcal{D}_k} 
\le \frac{2(K+K_1)}{k} \left\| \frac{R_k}{G_0^3} \right \|_{\infty, \mathcal{D}_k \cup \mathcal{E}_k }
+ \frac{2 K_8}{\beta_k} \left [\| \eta^{3/2} G_k \|_\infty 
+ \| \eta^{5/2} G_k^\prime \|_{\infty} \right ] $$
\end{Lemma}

\begin{proof} From (\ref{5.13}), since  $\eta_{j,k}$ are large 
 and therefore $g_j (\eta_{j',k})/g_j (\eta_{j,k})$ are exponentially
small  in $\beta_k$ for $j' \ne j$, it is clear that 
$$a_j g_j (\eta_{j,k}) \sim 
\mathcal{V} \left [\frac{R_k}{G_0^3} - b_{3,2} G_k^\prime - b_{3,1} G_k \right ] (\eta_{j,k})
- G_k (\eta_{j,k})$$
From Lemma (\ref{L5.2}), Corollaries \ref{L5.4} and \ref{L5.5}, it follows that
\begin{equation}
\label{5.17}
|a_j \eta_{j,k}^{3/2} g_j (\eta_{j,k})| <2 |\eta_{j,k}|^{3/2} |G_k (\eta_{j,k} )| + \frac{2 K_1}{\beta_k}
\Big\| \frac{R_k}{G_0^3} \Big\|_{\infty, \mathcal{D}_k} + 
\frac{{\hat K}_8}{\beta_k} \left [\| \eta^{3/2} G_k \|_{\infty,\mathcal{D}_k} 
+ \| \eta^{5/2} G_k^\prime \|_{\infty,\mathcal{D}_k} \right ] 
\end{equation}
Now, we conclude from Theorem~\ref{T4.1}  that
\begin{equation}
\label{115}
|\eta_{j,k}|^{3/2} |G_k (\eta_{j,k} )| \le \frac{K}{k} \Big\| \frac{R_k}{G_0^3} \Big\|_{\infty, \mathcal{E}_k}
\end{equation}
Since $\eta^{3/2} g_j (\eta)/(\eta_{j,k}^{3/2} g_j (\eta_{j,k}))$ are bounded  
independently of $\beta_k$ and the  proof  follows.
\end{proof}

\begin{Lemma}
\label{L5.9}
$$ \| \eta^{5/2} \mathcal{T}_2 \left [ G_k, G_k^\prime \right ] \|_{\infty, \mathcal{D}_k}
\le  \frac{{|\eta_{j,k}}^{3/2}|}{\beta_k^{2/3}}
\left \{ \frac{2(K+K_1)}{k} \Big\| \frac{R_k}{G_0^3} \Big\|_{\infty, \mathcal{D}_k \cup \mathcal{E}_k }
+ \frac{2 {\hat K}_8}{\beta_k} \left [\| \eta^{3/2} G_k \|_\infty 
+ \| \eta^{5/2} G_k^\prime \|_{\infty} \right ] \right \} $$
\end{Lemma}

\begin{proof} Taking into account the behavior of $h_j(\eta)$ for large $\eta$ we note that
$$ | \beta_k^{1/3} \eta_{j,k}^{5/2} h_j (\eta_{j,k}) a_j g_j (\eta_{j,k}) | \le  
C \beta_k^{1/3} |\eta_{j,k}|^{3/2}  | \eta_{j,k}^{3/2} a_j g_j (\eta_{j,k}) | $$
Using (\ref{5.17}) and (\ref{115}), the  proof  follows.
\end{proof}

\noindent{\bf Proof of Theorem \ref{T5.1}}. 
We consider the space Banach $\mathcal{B}$ 
of  pairs of analytic functions $(u, v)$ in the interior of 
$\mathcal{D}_k$ continuous in its closure with the  norm
$$
\| (u, v) \| = \| \eta^{3/2} u \|_{\infty, \mathcal{D}_k} + \| \eta^{5/2}
u^\prime \|_{\infty, \mathcal{D}_k} $$
Associating $G_k$ and $G_k^\prime$ in (\ref{5.10}) and (\ref{5.14}) with $u$
and $v$, 
we define the linear operator $L$ from $\mathcal{B}$ to $\mathcal{B}$
by
$$ L \left [ (u, v) \right ] 
= \left ( \mathcal{V} \left [{\hat R}_k [u,v] \right ] + 
\mathcal{T}_1 [u, v] , \mathcal{V'} \left [{\hat R}_k [u,v] \right ] 
+ \mathcal{T}_2 [u, v] \right ) $$ 
where ${\hat R}_k $ is now thought of as an
operator on $(u, v)$ for fixed $\frac{R_k}{G_0^3}$ such that
${\hat R}_k [G_k, G_k^\prime] (\eta) $ equals the right hand
side of (\ref{5.11}).  

It is a simple application of Lemmas 
\ref{L5.2}-\ref{L5.3}, \ref{L5.8}-\ref{L5.9} and
Corollaries \ref{L5.4} and \ref{L5.7} that 
$$ \| L \left [ (u, v) - ({\tilde u}, {\tilde v} \right) ]  \| 
\le {\tilde \delta} 
\| (u, v) - ({\tilde u},
{\tilde v} ) \| $$  
where 
$${\tilde \delta} = \max \left \{\frac{K_3}{\beta_k}, \frac{{\hat K}_8}{\beta_k} , 
\frac{K_4
\eta_{3,k}^{3/2}}{\beta_k^{2/3}}, \frac{{\hat K}_8 \eta_{3,k}^{3/2}}{
\beta_k^{2/3}}  \right \} <\frac{1}{2}$$
for sufficiently
large $\beta_k$ and small $\epsilon$.  
Contractivity of $L$ implies that it has a unique fixed point. 
The estimates in the Lemma follow from
(\ref{5.10}) and (\ref{5.14}). 

\noindent{\bf Proof of Lemma \ref{L0.2}}. 
First, for $k=1,...,k_0$, the statement in the Lemma holds if $A$
is sufficiently large (depending on $k_0$) in a common domain
$\mathcal{D}_{k_0} \cup \mathcal{E}_{k_0} $, chosen to contain
$\mathcal{D}_{k_0+1} \cup \mathcal{E}_{k_0+1}$. Assume therefore
  that $k >k_0$ where $k_0+1$ is large enough to ensure contractivity
  in Theorems \ref{T4.1} and \ref{T5.1}.  Assume the  statement 
  holds $j=1,...,k_0$ in a common domain $\mathcal{D}_{k_0} \cup
  \mathcal{E}_{k_0} $ and for $j=k_0+1$, ..., $k-1$ in a corresponding
  sequence of domains $\mathcal{D}_j \cup \mathcal{E}_j$. It follows
  from the construction of these domains that it then holds in
  $\mathcal{D}_k \cup \mathcal{E}_k$.  We then get the estimates on
  $R_k$ needed in Theorems \ref{T4.1} and \ref{T5.1}, which imply
$$ \| \eta^{3/2} G_k \|_{\infty, \mathcal{D}_k \cup \mathcal{E}_k}  
\le \frac{K_{10}}{k} \Big\| \frac{R_k} {G_0^3} \Big\|_{\infty, \mathcal{D}_k 
\cup \mathcal{E}_k } $$
$$ \| \eta^{5/2} G_k^\prime \|_{\infty, \mathcal{D}_k \cup \mathcal{E}_k}  
\le K_{11} \Big\| \frac{R_k} {G_0^3} \Big\|_{\infty, \mathcal{D}_k \cup \mathcal{E}_k} 
$$
and therefore, from the estimates on $\| R_k \eta^{3/2} \|$ in (\ref{0.3.7}), we get
$$ \| \eta^{3/2} G_k \|_{\infty, \mathcal{D}_k \cup \mathcal{E}_k}  
\le \frac{K_{10} K_3}{k^3} (B^2 A^k + B A^{k-1})  $$ 
$$ \| \eta^{5/2} G_k^\prime \|_{\infty, \mathcal{D}_k \cup \mathcal{E}_k}  
\le \frac{K_{11} K_3 }{k^2} (B^2 A^k + B A^{k-1} ) $$ 
Using   eq. (\ref{0.2.8}) and the bounds on
$R_k $, it follows that 
\begin{multline*}
\| G_k^{\prime \prime \prime} \|_{\infty, \mathcal{D}_k \cup \mathcal{E}_k }
\le 
K_{12} \| \eta^{3/2} G_k \|_{\infty, \mathcal{D}_k \cup \mathcal{E}_k}
+ K_{13} \| \eta^{5/2} G_k^\prime \|_{\infty, \mathcal{D}_k \cup \mathcal{E}_k}
+ K_{14} \| \eta^{3/2} R_k \|_{\infty, \mathcal{D}_k \cup \mathcal{E}_k} \\
\le \left ( K_{12} K_{11} + \frac{K_{13}}{k} K_{10} + K_{14} \right ) 
\frac{K_3}{k^2} (B^2 A^k + B A^{k-1}) 
\end{multline*} 
It is clear that for $B$ sufficiently small and $A$ sufficiently large,  
the estimates
(\ref{0.3.4})-(\ref{0.3.6}) on $G_k$, $G_k^\prime$ and $G_k^{\prime \prime
  \prime}$  follow. 
The result follows now by induction.

\noindent{\bf Proof of Theorem \ref{T0.1}}.  Now  this follows easily
from Lemma \ref{L0.2} since the estimates guarantee convergence of the Taylor
series (\ref{0.2.6}) for sufficiently small $\tau$.

\section{Appendix: Singularities of nonlinear ODEs}

We first mention briefly a number of results in \cite{Invent} and then allow
for slight modifications in the assumptions, to adjust  for the equation of
$G_0$.

\subsection{Setting of \cite{Invent} and generalizations}\label{Sett} We adopt, with few exceptions that we mention, the same conditions,
notations and terminology as \cite{DMJ} and \cite{Invent}; the results on
formal solutions and their generalized Borel summability are also taken from
\cite{DMJ}.

The differential system considered has the form

\begin{eqnarray}
 \label{eqor}
  \mathbf{y}'=\mathbf{f}(x^{-1},\mathbf{y})  \qquad
  \mathbf{y}\in\CC^n,\ \ x\in\CC              
   \end{eqnarray}
   
   \z where 
   
 \z   (i) $\mathbf{f}$ is {\em analytic} in a neighborhood
   $\mathcal{V}_x\times\mathcal{V}_\mathbf{y}$ of $ (0,\mathbf{0})$, under the
   genericity conditions that:
   
   \z (ii) the eigenvalues $\lambda_j$ of the matrix
   $\hat\Lambda=-\left\{\frac{\partial f_i}{\partial
       y_j}(0,\mathbf{0})\right\}_{i,j=1,2,\ldots n}$ are linearly independent
   over $\ZZ$ (in particular $\lambda_j\ne 0$) and such that $\arg\lambda_j$
   are all different.

 We now allow for the same assumptions, except we replace (ii) by

  \z (ii') There is at most one zero eigenvalue of $\hat\Lambda$ and all the
  other $\lambda_j$ are linearly independent over $\ZZ$ (in particular
  $\lambda_j\ne 0$) and such that $\arg\lambda_j$ are all different.

By elementary changes of variables, the system (\ref{eqor}) can be
brought to the {\em normalized form} \cite{DMJ}.

\begin{eqnarray}\label{eqor1}
{\bf y}'=-\hat\Lambda {\bf y}+
\frac{1}{x}\hat A {\bf y}+{\bf g}(x^{-1},{\bf y})
\end{eqnarray}

\z where $\hat{\Lambda}=\mbox{diag}\{\lambda_j\},\ 
\hat{A}=\mbox{diag}\{\alpha_j\}$ are constant matrices, $\mathbf{g}$ is
analytic at $(0,\mathbf{0})$ and ${\bf g}(x^{-1},{\bf y})=
O(x^{-2})+O(|\bfy|^2)$ as $x\rightarrow\infty$ and
$\mathbf{y}\rightarrow 0$. 

As in \cite{DMJ} we normalize the system so that $\Re(\alpha_j)>0$.

Performing a further transformation of the
type $\bfy\mapsto \bfy -\sum_{k=1}^{M}\mathbf{a}_k x^{-k}$ (which takes
out $M$ terms of the formal asymptotic series solutions of the
equation), makes

\begin{equation}{\label{eqggg}}{\bf g}(|x|^{-1},{\bf y})= O(x^{-M-1};|\bfy|^2;|x^{-2}\bfy|)\ \ \ \ 
(x\rightarrow\infty;\ \bfy\rightarrow 0)\end{equation}
where
$$M\ge\max_j\Re(\alpha_j)$$
and $O(a;b;c)$ means (at most) of the order of the largest
among $a,b,c$.

Our analysis applies to solutions ${\mathbf{y}(x)}$ such that
$\mathbf{y}(x)\rightarrow 0$ as $x\rightarrow\infty$ along some
arbitrary direction $d=\{x\in\CC:\arg(x)=\phi\}$. A movable singularity of
$\mathbf{y}(x)$ is a point $x\in\CC$ with $x^{-1}\in\mathcal{V}_x$
where $\mathbf{y}(x)$ is not analytic. The point at infinity is an irregular
singular point of rank 1; it is a fixed singular point of the system
since, after the substitution $x=z^{-1}$ the r.h.s of the transformed
system, $\frac{dy}{dz}=-z^{-2}\mathbf{f}(z,\mathbf{y})$ has, under the
given assumptions, a pole at $z=0$.

An $n$-parameter formal solution of (\ref{eqor1}) (under the assumptions
mentioned) as a combination of powers and exponentials is found in the form
\begin{gather}
  \label{transs}
  \tilde{\mathbf{y}}(x)=\sum_{\mathbf{k}\in (\NN\cup\{0\})^n}
  \mathbf{C}^{\mathbf{k}}\erm^{-\boldsymbol{\lambda}\cdot\mathbf{k}x}
  x^{\boldsymbol{\alpha}\cdot\mathbf{k}}\tilde{\mathbf{s}}_{\mathbf{k}}(x)
\end{gather}

\z where $\tilde{\mathbf{s}}_\mathbf{k}$ are (usually factorially divergent)
formal power series: $\tilde{ \mathbf{s}}_{\bf{0}}=\tilde{
  \mathbf{y}}_{\bf{0}}$ and in general
\begin{gather}\label{expds}
  \tilde{ \mathbf{s}}_\mathbf{k}(x)=\sum_{r=0}^{\infty}\frac{\tilde{
      \mathbf{y}}_{\mathbf{k};r}}{x^{r}}
\end{gather}
that can be determined by formal substitution of (\ref{transs}) in
(\ref{eqor1}); $\mathbf{C}\in\CC^n$ is a vector of
parameters\footnote{In the general case when some assumptions made
  here do not hold, the general formal solution may additionally
logs  iterated exponentials,  and powers \cite{Ecalle}. The
  present paper only discusses equations in the setting explained at
  the beginning of the present section.}(we
use the notations $\mathbf{C}^{\mathbf{k}}=\prod_{j=1}^n C_j^{k_j}$,
$\boldsymbol{\lambda}=(\lambda_1,...,\lambda_n)$,
$\boldsymbol{\alpha}=(\alpha_1,...,\alpha_n)$,
$|\mathbf{k}|=k_1+...+k_n$).

Note the structure of (\ref{transs}): an infinite sum of (generically)
divergent series multiplying exponentials. They are called {\em{formal
    exponential power series}} \cite{Wasow}.

From the point of view of correspondence of these formal solutions to
actual solutions it was recognized that not all expansions
(\ref{transs}) should be considered meaningful; also they 
are defined relative to a sector (or a direction).  

Given a direction $d$ in the complex $x$-plane the {\em{transseries}}
(on $d$), introduced by \'Ecalle \cite{Ecalle}, are, in our context,
those exponential series (\ref{transs}) which are formally {\em
  asymptotic} on $d$, i.e. the terms
$\mathbf{C}^{\mathbf{k}}\erm^{-\boldsymbol{\lambda}\cdot\mathbf{k}x}
x^{\boldsymbol{\alpha}\cdot\mathbf{k}}x^{-r}$ (with $\mathbf{k}\in
(\NN\cup\{0\})^n,\, r\in\NN\cup\{0\}$) form a well ordered set with
respect to $\gg$ on $d$ (see also \cite{DMJ}).\footnote{We note here a
  slight difference between our transseries and those of \'Ecalle, in
  that we are allowing complex constants.} (For example, this is the
case when the terms of the formal expansion become (much) smaller when
$\mathbf{k}$ becomes larger.)

We recall that the {\em{antistokes lines}} of (\ref{eqor1}) are the
$2n$ directions of the $x$-plane $i\overline{\lambda_j}\,\RR_+,\ -
i\overline{\lambda_j}\, \RR_+,\ j=1,...,n$, i.e. the directions along
which some exponential $\e^{-\lambda_jx}$ of the general formal
solution (\ref{transs}) is purely oscillatory.

In the context of differential systems with an irregular singular
point, asymptoticity should be (generically) discussed relative to a
direction towards the singular point; in fact, under the present
assumptions (of non-degeneracy) asymptoticity can be defined on
sectors.

Let $d$ be a direction in the $x$-plane which is not an antistokes
line. The solutions
$\mathbf{y}(x)$ of (\ref{eqor1}) which satisfy
\begin{eqnarray}
  \label{eq:defasy0}
 \mathbf{y}(x)\rightarrow 0\ \  (x\in d;\ |x|\rightarrow\infty)
\end{eqnarray}

\z are analytic for large $x$ in a sector containing $d$, between two
neighboring antistokes lines and have the same asymptotic series

\begin{eqnarray}
  \label{eq:asy0}
  \mathbf{y}(x)\sim \tilde{\mathbf{y}}_{\bf{0}}\ \  (x\in d;\ |x|\rightarrow\infty)
\end{eqnarray}

In the context of (\ref{eqor1}), a generalized Borel summation
$\mathcal{L}\mathcal{B}$ of transseries (\ref{transs}) is defined in \cite{DMJ}.

The formal solutions (\ref{transs}) are determined by the equation
(\ref{eqor1}) that they satisfy, except for the parameters $\bf
C$. Then a correspondence between actual and formal solutions of the
equation is an association between solutions and constants $\bf
C$. This is done using a generalized Borel summation $\mathcal{L}\mathcal{B}$.

The operator $\mathcal{L}\mathcal{B}$ constructed in \cite{DMJ} can be
applied to any transseries solution (\ref{transs}) of (\ref{eqor1})
(valid on its open sector $S_{trans}$, assumed non-empty) on
any direction $d\subset S_{trans}$ and yields an actual solution
$\mathbf{y}=\mathcal{L}\mathcal{B}\tilde{\mathbf{y}}$ of
(\ref{eqor1}), analytic in a domain $S_{an}$.
Conversely, any solution ${\mathbf{y}}(x)$ satisfying (\ref{eq:asy0})
on a direction $d$ is represented as
$\mathcal{LB}\tilde{\mathbf{y}}(x)$, on $d$, for some unique
$\tilde{\mathbf{y}}(x)$:

\begin{equation}
 \label{transsf}
  \mathbf{y}(x)=
\sum_{\mathbf{k}\ge 0}
  \mathbf{C}^{\mathbf{k}}\erm^{-\boldsymbol{\lambda}\cdot\mathbf{k}x}
  x^{\mathbf{M}\cdot\mathbf{k}}\mathbf{y}_{\mathbf{k}}(x) 
=\sum_{\mathbf{k}\ge 0}
  \mathbf{C}^{\mathbf{k}}\erm^{-\boldsymbol{\lambda}\cdot\mathbf{k}x}
  x^{\mathbf{M}\cdot\mathbf{k}}\lap\bor\tilde{\mathbf{y}}_{\mathbf{k}}(x)=\mathcal{L}\mathcal{B}\tilde{\mathbf{y}}(x)
\end{equation}
for some constants $\mathbf{C}\in\CC^n$, where $M_j=\lfloor
\Re\alpha_j\rfloor+1$ ($\lfloor \cdot \rfloor$ is the integer part), and
\begin{equation}\label{expd}
  \tilde{\mathbf{y}}_\mathbf{k}(x)=\sum_{r=0}^{\infty}\frac{\tilde{\mathbf{y}}_{\mathbf{k};r}}{x^{-\mathbf{k}\boldsymbol{\alpha}'+r}}\ \ \ \ \ \ \ \
 (\boldsymbol{\alpha}'=\boldsymbol{\alpha}-\bf M)
\end{equation}
(for technical reasons the Borel summation procedure is applied to the series
\begin{equation}\label{relyksk}
\tilde{ \mathbf{y}}_\mathbf{k}(x)= x^{\mathbf{k}\boldsymbol\alpha'}
  \tilde{ \mathbf{s}}_\mathbf{k}(x)
\end{equation}
\z rather than to $ \tilde{ \mathbf{s}}_\mathbf{k}(x)$ cf.
(\ref{transs}),(\ref{expds})).

The modification necessary to extend (\ref{transsf}) to the case $\lambda_0=0$
is outlined in \S\ref{ext}.
\subsection{Normal form of Eq. (\ref{0.2.7})}

We first give some detail on the normalization procedure, in the limit
$|x|\rightarrow\infty$. It can be checked that there is a one-parameter family
of formal solutions to (\ref{0.2.7}) in the form $Cx^{-1/2}-\frac{15
  C^4}{8}x^{-5}+...$. The physical problem requires $C=1$; this suggests the
substitution $G_0=x^{-1/2}+h(x)$ where $h$ is expected to behave like
$-\frac{15}{8}x^{-5}$. 

The normalizing substitution produces an equation with solutions in the form
(\ref{transsf}), where the terms with $\bfk>0$ contain exponentials with
argument linear in the final variable; the type of the exponenential in the
equation for $h$ can be found by linear perturbation theory around a solution
$h_0$; with $h-h_0=\delta$, the leading order equation for $\delta$ is

\begin{equation}
  \label{edelta}
  \delta'''+\frac{2}{9}x^{5/2}\delta'+\frac{1}{9}x^{2/3}\delta=0
\end{equation}
where the substitution of the form $\delta=A(x)\e^{bx^p}$ shows that $p=9/4$
implying that the natural variable is $x^{9/4}$.

Taking $G_0=x^{-1/2}+x^{-1/2}g(x^{9/4})$, $\xi=x^{9/4}$ in (\ref{0.2.7})
we obtain
\begin{equation}
  \label{eg3}
  g'''+\frac{1}{\xi}g''+\left(\frac{11}{81\xi^2}+\frac{32}{729}\frac{1}{(1+g)^3}\right)g'=\frac{40}{243}\left(\frac{1}{\xi^3}+\frac{g}{\xi^3}\right)
\end{equation}

\z which, written as a system, becomes

\begin{equation}
  \label{matrf}
{ \begin{pmatrix}
g''\\g'\\g
\end{pmatrix}}'=
  \begin{pmatrix}
0&-\frac{32}{729}&0\\1&0&0\\0&1&0
\end{pmatrix}{ \begin{pmatrix}
g''\\g'\\g
\end{pmatrix}}-\frac{1}{\xi} \begin{pmatrix}
1&0&0\\0&0&0\\0&0&0
\end{pmatrix}{ \begin{pmatrix}
g''\\g'\\g
\end{pmatrix}}+O(g^2,\xi^{-2})
\end{equation}

\z The eigenvalues of the first matrix on the rhs of (\ref{matrf}),
$\{0,\pm\frac{4i\sqrt{2}}{27}\}$, are the values of $\lambda$ in (\ref{transsf}).
The fact that one eigenvalue is zero requires a slight modification
in the proofs of \cite{DMJ}.

\subsection{Extension of the proofs in \cite{DMJ} to the assumption (ii')}\label{ext}

In an attempt to minimize the possibility of confusion with the setting
in \cite{DMJ} we assume  that  the order of the system is $n+1$, 
we  count dimensions  starting with zero, and take $\lambda_0=0$. There
is no contribution from $\lambda_0$ to the general formally decreasing
transseries (\ref{transs}); this is due to the normalization
$\Re(\alpha_j)>0$.
  
 The convolution equations satisfied by $\mathbf{Y}=\bor\mathbf{y}$ and
  $\mathbf{Y}_{\bf k}=\bor\mathbf{y}_{\bf k}$ are given
  still given by equations (1.13 ) and (1.16) as in \cite{DMJ} (with the
  notation
$\hat{A}=-\hat{B}$ used there):

\begin{eqnarray}\label{eqil}
-p{\bf Y}={\bf F}_0-\hat\Lambda {\bf Y}-
\hat B\mathcal{P}\mathbf{Y}+{\cal  N}({\bf Y})
\end{eqnarray}

\begin{eqnarray}
  \label{invlapvk}\label{eqMv}\label{invlapyk}
  &&\left(-p+\hat{\Lambda}-\bfk\cdot\boldsymbol{\lambda}\right)\bfY_\bfk
+\left(\hat{B}+\bfk\cdot\bfm\right)
\mathcal{P}\bfY_\bfk+\sum_{|\bfj|=1}\bfd_\bfj*\bfY_\bfk^{*\bfj}
=\bfT_\bfk\cr&&
\end{eqnarray}

The only difference relevant to \cite{Invent} with respect to the analysis in
\cite{DMJ} is in the study of $\bfY$, and once the analog results are
obtained, the analysis of $\bfY_\bfk$ is virtually identical. By the
normalization choice, we have $\mathbf{F_0}=p^{M}\mathbf{H}(p)$ where
$\mathbf{H}$ is analytic at zero.  In the equation (2.35) of \cite{DMJ}

\begin{eqnarray}\label{eqilm}
\bfY=
\left(\hat{\Lambda}-p\right)^{-1}\left({\bf F}_0-
\hat B\mathcal{P}\bfY+{\cal  N}({\bf Y})\right)=\mathcal{M}(\bfY)
\end{eqnarray}

\z we separate the zeroth component which is apparently singular (as was done
in the study of $\mathbf{Y}_1$ in\cite{DMJ} \S 2.2.2; here the analysis is simpler):

\begin{eqnarray}\label{eqilm1}
-p(\bfY)_0-\alpha_0\mathcal{P}(\bfY)_0
={\bf F}_{0;0}+\left({\cal  N}({\bf Y}) \right)_0:=\mathbf{R}_0
\end{eqnarray}

\z or

\begin{eqnarray}\label{eqilm1}
-p(\bfY)'_0-(\bfY)_0-\alpha_0(\bfY)_0=\mathbf{R}_0'
\end{eqnarray}

\z which we rewrite as an integral equation, which after integration by parts
reads:

\begin{eqnarray}\label{eqilm1}
(\bfY)_0=-{\bf F}_{0;0}+(1+\alpha_0)\int_0^1{\bf F}_{0;0}(tp)dt-\left({\cal  N}({\bf Y}) \right)_0+(1+\alpha_0)\int_0^1\left({\cal  N}({\bf Y}) \right)_0(tp)dt=\mathcal{M}_0^{[1]}(\mathbf{Y})
\end{eqnarray}

\z The system is of the form (\ref{eqilm})

\begin{eqnarray}\label{eqilm3}
\bfY=
\mathcal{M}^{[1]}(\bfY)
\end{eqnarray}
with $\mathcal{M}^{[1]}=\mathcal{M}$ for all components other than the zeroth
one defined in (\ref{eqilm1}). The equation (\ref{eqilm3}) is contractive in
the ball $B=\{\bfY:\{p:|p|<\epsilon\}:\|\bfY\|_\infty<2\epsilon\}$ for small
enough $\epsilon$, and also in the focusing algebra (3a) in \S2.1.1 in
\cite{DMJ} for $\beta_k=1$ (allowed by the normalization of $\mathbf{F}_0$) as
follows from immediate estimates.

No other nontrivial adaptations are needed in the proofs in \cite{DMJ}.

\subsection{Results of \cite{Invent} as extended in \S\ref{ext}}

The map $\tilde{\mathbf{y}}\mapsto\mathcal{L}\mathcal{B}(\tilde{\mathbf{y}})$
depends on the direction $d$, and (typically) is discontinuous at the finitely
many Stokes lines, see \cite{DMJ}, Theorem 4.

For linear equations only the directions $\overline{\lambda_j}\,\RR_+,\ 
j=1,...,n$ are Stokes lines, but for nonlinear equations there are also other
Stokes lines, recognized first by \'Ecalle.  $\mathcal{LB}$ is only
discontinuous because of the jump discontinuity of the vector of ``constants''
$\mathbf{C}$ across Stokes directions (Stokes' phenomenon); between Stokes
lines $\mathcal{L}\mathcal{B}$ does not vary with $d$.

The function series in (\ref{transsf}) is uniformly {\em convergent}
and the functions $\mathbf{y}_\bfk$ are analytic on domains $S_{an}$ (for some $\delta>0$,
$R=R(\mathbf{y}(x),\delta)>0$.

\begin{Theorem}
\label{T1}
There exists $\delta_1 >0$ so that for $|\xi|<\delta_1$ the
power series

\begin{equation}
  \label{eq:formFm}
  \mathbf{F}_m(\xi)=\sum_{k=0}^{\infty}\xi^{k}\tilde{\mathbf{y}}_{k\mathbf
  {e}_1;m},\ \ m=0,1,2,...
\end{equation}

\z converge. Furthermore

 \begin{gather}
   \label{estiy}
   \mathbf{y}(x)\sim
   \sum_{m=0}^{\infty}x^{-m}\mathbf{F}_m(\xi(x)) \ \ 
   (x\in\mathcal{S}_{\delta_1}, \ x\rightarrow\infty)
 \end{gather}

 \z uniformly in $\mathcal{S}_{\delta_1}$, and  the asymptotic representation
 (\ref{estiy}) is differentiable.
 
 The functions $\mathbf{F}_m$ are uniquely defined by (\ref{estiy}),
 the requirement of analyticity at $\xi=0$, and
 $\mathbf{F}_0'(0)=\mathbf{e}_1$.

\end{Theorem}

\begin{Remark}\label{R=1} A direct calculation shows that the functions
  $\mathbf{F}_m$ are solutions of the system of equations

%remarkR=1

\begin{align}\label{eqF0princ}
  &\frac{\mathrm{d}}{\mathrm{d}\xi}\mathbf{F}_0=\xi^{-1}\left(\hat{\Lambda}\mathbf{F}_0-\mathbf{g}(0,\mathbf{F}_0)\right)
  \\ \label{eq:system1}
  &\frac{\mathrm{d}}{\mathrm{d}\xi}\mathbf{F}_m+\hat{N}\mathbf{F}_m=
 \alpha_1\frac{\mathrm{d}}{\mathrm{d}\xi}
  \mathbf{F}_{m-1}+\mathbf{R}_{m-1}
\ \ \  \ \  {\mbox{for}}\ m\ge 1
\end{align}

\z where $\hat{N}$ is the matrix 

\begin{equation}\label{defN}
\xi^{-1}(\partial_{\mathbf{y}}
\mathbf{g}(0,\mathbf{F}_0)-\hat{\Lambda})
\end{equation}
\z and the function $\mathbf{R}_{m-1}(\xi)$ depends only on the
$\mathbf{F}_k$ with $k<m$:

 \begin{align}\label{R**}
  \left.\xi\mathbf{R}_{m-1}=
  -\left[ (m-1)I+\hat{A}\right] \mathbf{F}_{m-1}
  -
\frac{1}{ m!}\frac{\mathrm{d}^m}{\mathrm{d}z^m}\mathbf{g}
  \left( z; \sum_{j=0}^{m-1} z^j \mathbf{F}_j \right) \right|_{z=0}
\end{align}

\end{Remark}

\subsection{Formal arguments for thin-film equation}

Consider the particular initial value problem in one space dimension:
\begin{equation}
\label{film2}
h_t + \left ( h^3 h_{xxx} \right )_{x} = 0 ~~~,~~h(x, 0) = \frac{1}{1+x^2}
\end{equation}
This is a special case of $h_t + (h^n h_{xxx})_x=0$.  Global existence
proofs are available only for $n > 3.5$; numerical solutions suggest
finite-time singularity for $n=1$ \cite{Bertozzi}.

For the problem (\ref{film2}) and variations of it, the complex region
for which existence is expected, at least for small $t$, includes the
real $x$-axis.  For the specific initial value problem, we change
variables:
$$
h(x, t) = H(1+x^2, t) ,~~\\~~ ~\xi = 1+x^2 $$
and obtain a nonlinear
PDE for $H(\xi, t)$. A formal asymptotic expansion in powers of $t$
results in
\begin{equation}
\label{eqformal}
H(\xi, t) = \frac{1}{\xi} \sum_{j=0}^\infty
P_{2j} \left (t^{1/2} \xi^{-7/2}, t^{1/2} \xi^{-5/2} \right )
\end{equation}
where $P_{2 j}$ are homogenous polynomials of order $2j$.  With
appropriate changes of variables, we expect the regularity theorem
\cite{invent03} to be adaptable to prove short term existence for a
complex $\xi$ sector that includes $(1, \infty)$ (i.e., $x\in\RR$),
and show further the validity of (\ref{eqformal}) for $\xi\gg t^{1/7}
$ in this sector.

Asymptotics (\ref{eqformal})
fails when $\xi = O(t^{1/7})$.
Introducing scaled variables,
$$ \eta = \xi t^{-1/7} , \tau = t^{1/7} \\,
H(\xi (\eta,\tau), t(\tau)) =
\xi^{-1} F(\eta,\tau),  $$
gives a formal solution as an expansion in integer powers of $\tau$,
\begin{equation}
\label{eqFseries}
F(\eta, \tau) =\sum_{k=0}^\infty \tau^k F_k (\eta)
\end{equation}
We expect this series to be convergent.
The equation of $F_0$ can be integrated once by
using far-field matching condition to give:
$$
F_{0}^3 F_0^{\prime\prime\prime} -\frac {6}{\eta^3} F_{0}^4 -\frac
{\eta^4}{112} F_0 + \frac{6}{\eta^2} F_0^3 F_0^\prime - \frac
{3}{\eta} F_{0}^3 F_0^{\prime \prime} +\frac {\eta^4}{112} = 0 $$
With
the further transformation $F_0 (\eta) = 1 +y(\eta^{7/3})$, the
equation for $y$ is in a form to which the general theory
\cite{Invent} applies. From the leading order singularity of the ODE,
and the expected convergence of (\ref{eqFseries}), as for modified
Harry-Dym, we expect to show that the thin-film equation has
singularities at points close to $x_s (t)$ with
$1+x_s^2 = \eta_s
t^{1/7}$.

\vfill \eject

\newpage
\figure 

\scalebox{0.5}[0.5]{
\includegraphics{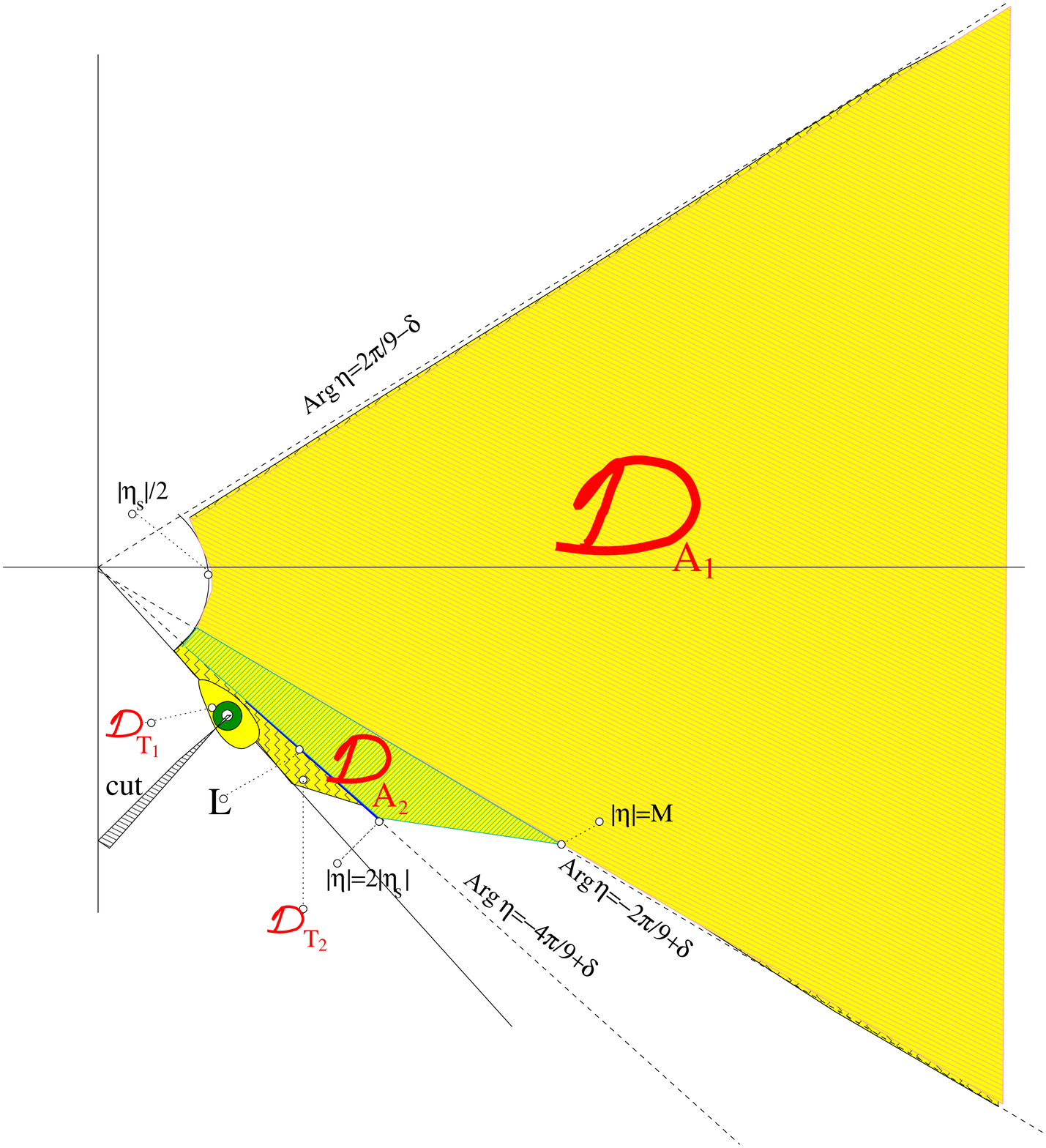}}
\caption{Subdomains of $\mathcal{D}$}
\endfigure

\newpage
\figure
\scalebox{0.5}[0.5]{
\includegraphics{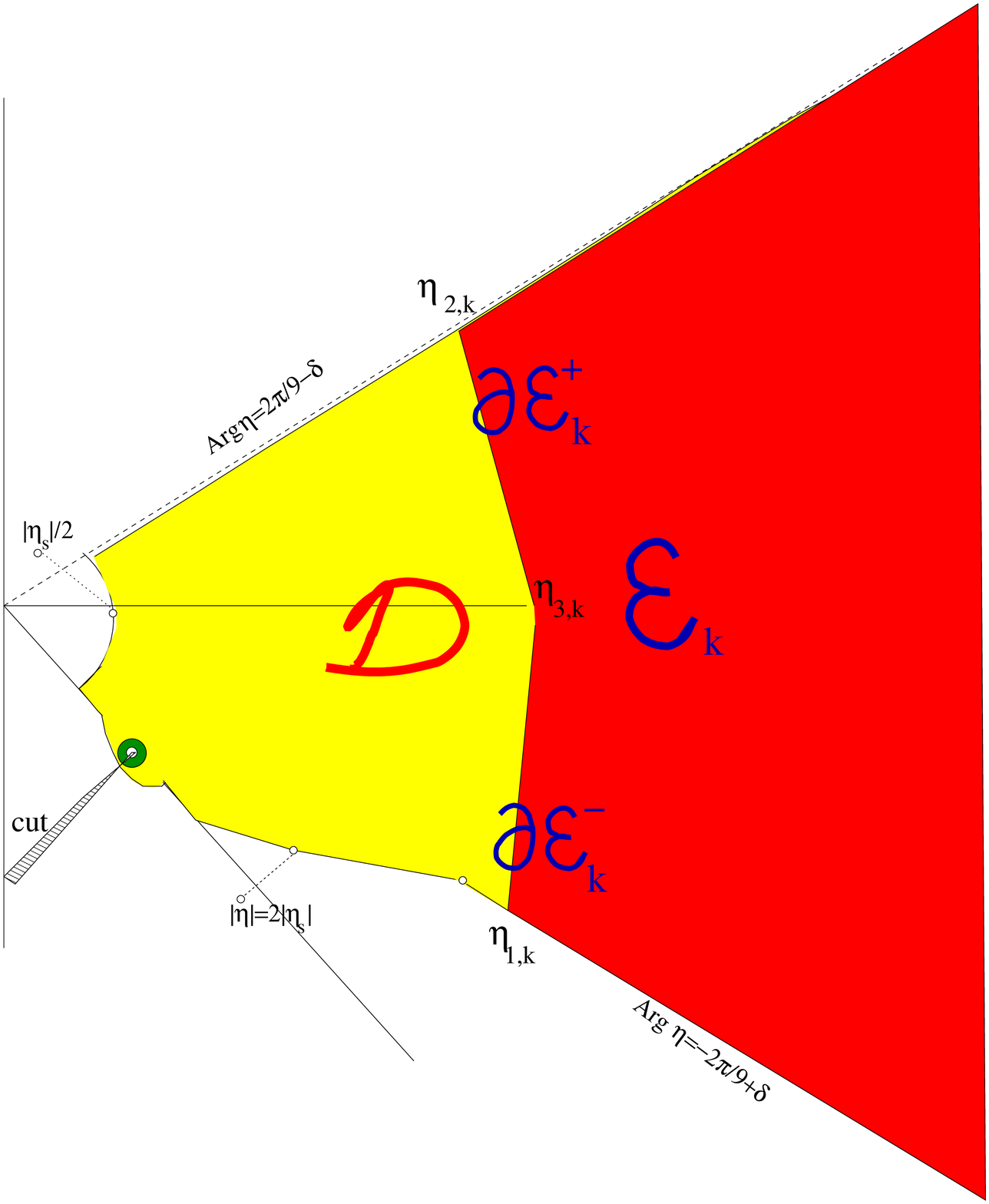}}
\caption{Domains $\mathcal{D}_k$, $\mathcal{E}_k$ and common boundary
$\partial{E}_k$.}
\endfigure

\newpage
\figure

\scalebox{0.5}[0.5]{
\includegraphics{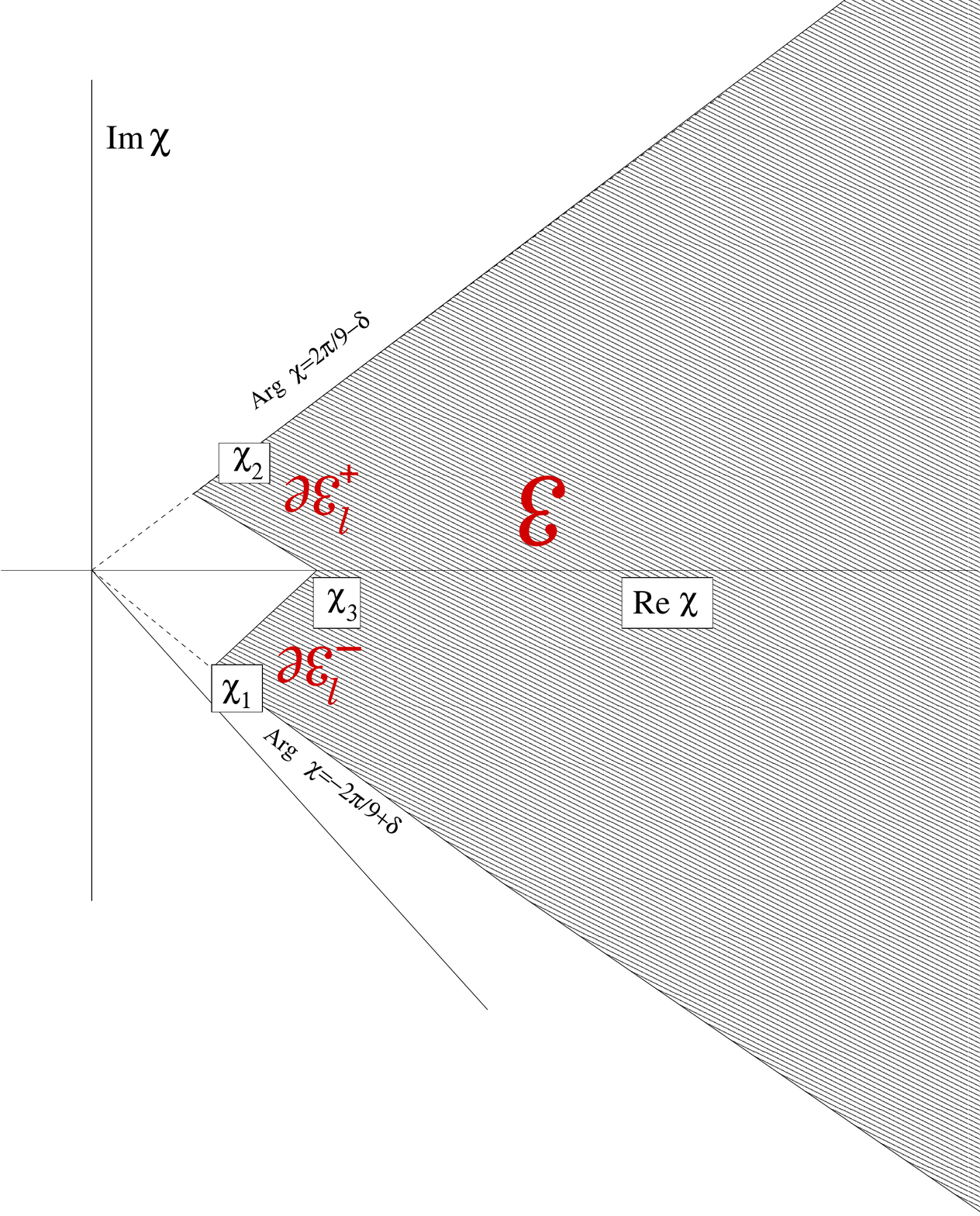}}
\caption{Domain $\mathcal{E}$ in the $\chi$-plane}
\endfigure

\end{document}